\documentclass[aap,MSNbibl,seceqn,citesort,dvips]{arximspdf}
\usepackage{subenv}

\doi{10.1214/09-AAP651}
\volume{20}
\issue{4}
\pubyear{2010}
\firstpage{1359}
\lastpage{1388}

\makeatletter

\newproclaim{definition}{Definition}[section]

\newtheorem{proposition}[definition]{Proposition}
\newtheorem{theorem}[definition]{Theorem}
\newtheorem{corollary}[definition]{Corollary}

\newproclaim{remark}[definition]{Remark}

\makeatother

\begin{document}
\begin{frontmatter}

\title{The emergence of rational behavior in the presence of
stochastic perturbations}
\runtitle{Rationality and stochastic perturbations}

\begin{aug}
\author[A]{\fnms{Panayotis} \snm{Mertikopoulos}\corref{}\thanksref{t1,t2}\ead[label=e1]{pmertik@phys.uoa.gr}} and
\author[A]{\fnms{Aris L.} \snm{Moustakas}\thanksref{t1}\ead[label=e2]{arislm@phys.uoa.gr}}
\runauthor{P. Mertikopoulos and A. L. Moustakas}
\affiliation{University of Athens}
\address[A]{Department of Physics\\
University of Athens\\
Panepistimioupolis, Zografou\\
Athens, 15784\\
Greece\\
\printead{e1}\\
\phantom{E-mail: }\printead*{e2}} 
\end{aug}

\thankstext{t1}{Supported in part by the European Commission
Grants EU-FET-FP6-IST-034413
(Net-ReFound) and EU-IST-NoE-FP6-2007-216715 (NewCom$++$) and the Greek
Research Council project \textit{Kapodistrias} (no. 70/3/8831).}

\thankstext{t2}{Supported by the Empirikeion Foundation of Athens, Greece.}

\received{\smonth{6} \syear{2009}}
\revised{\smonth{9} \syear{2009}}

%
\begin{abstract}
We study repeated games where players use an exponential learning
scheme in order to adapt to an ever-changing environment. If the game's
payoffs are subject to random perturbations, this scheme leads to a new
stochastic version of the replicator dynamics that is quite different
from the ``aggregate shocks'' approach of evolutionary game theory.
Irrespective of the perturbations' magnitude, we find that strategies
which are dominated (even iteratively) eventually become extinct and
that the game's strict Nash equilibria are stochastically
asymptotically stable. We complement our analysis by illustrating these
results in the case of congestion games.
\end{abstract}

%

\end{frontmatter}

\section{Introduction}
\label{sec:introduction}

Ever since it was introduced in \cite{Na51}, the notion of a Nash
equilibrium and its refinements have remained among the most prominent
solution concepts of noncooperative game theory. In its turn, not only
has noncooperative game theory found applications in such diverse
topics as economics, biology and network design, but it has also become
the standard language to actually \textit{describe} complex agent
interactions in these fields.

Still, the issue of why and how players may arrive to equilibrial
strategies in the first place remains an actively debated question.
After all, the complexity of most games increases exponentially with
the number of players and, hence, identifying a game's equilibria
quickly becomes prohibitively difficult. Accordingly, as was first
pointed out by Aumann in \cite{Au74}, a player has no incentive to play
his component of a Nash equilibrium unless he is convinced that all
other players will play theirs. And if the game in question has
multiple Nash equilibria, this argument gains additional momentum: in
that case, even players with unbounded deductive capabilities will be
hard-pressed to choose a strategy.

From this point of view, rational individuals would appear to be more
in tune with Aumann's notion of a correlated equilibrium where
subjective beliefs are also taken into account \cite{Au74}.
Nevertheless, the seminal work of Maynard Smith on animal conflicts
\cite{MS74} has cast Nash equilibria in a different light because it
unearthed a profound connection between evolution and rationality:
roughly speaking, one leads to the other. So, when different species
contend for the limited resources of their habitat, evolution and
natural selection steer the ensuing conflict to an equilibrial state
which leaves no room for irrational behavior. As a consequence,
instinctive ``fight or flight'' responses that are deeply ingrained in
a species can be seen as a form of rational behavior, acquired over the
species' evolutionary course.

Of course, this evolutionary approach concerns large populations of
different species which are rarely encountered outside the realm of
population biology. However, the situation is not much different in the
case of a finite number of players who try to learn the game by playing
again and again and who strive to do better with the help of some
learning algorithm. Therein, evolution does not occur as part of a
birth/death process; rather, it is a byproduct of the players' acquired
experience in playing the game---see \cite{FL98} for a most
comprehensive account.

It is also worth keeping in the back of our mind that in some
applications of game theory, ``rationality'' requirements precede
evolution. For example, recent applications to network design start
from a set of performance aspirations (such as robustness and
efficiency) that the players (network devices) seek to attain in the
network's equilibrial state. Thus, to meet these requirements, one has
to literally reverse-engineer the process by finding the appropriate
game whose equilibria will satisfy the players---the parallel with
mechanism design being obvious.

In all these approaches, a fundamental selection mechanism is that of
the \textit{replicator dynamics} put forth in \cite{TJ78} and \cite{SS83}
which reinforces a strategy proportionately to the difference of its
payoff from the mean (taken over the species or the player's
strategies, depending on the approach). As was shown in the
multi-population setting of Samuelson and Zhang \cite{SZ92} (which is
closer to learning than the self-interacting single-population scenaria
of \cite{TJ78} and \cite{SS83}), these dynamics are particularly
conducive to rationality. Strategies that are suboptimal when paired
against any choice of one's adversaries rapidly become extinct, and in
the long run, only rationally admissible strategies can survive. Even
more to the point, the only attracting states of the dynamics turn out
to be precisely the (strict) Nash equilibria of the game---see \cite
{We95} for a masterful survey.

We thus see that Nash equilibria arise over time as natural attractors
for rational individuals, a fact which further justifies their
prominence among noncooperative solution concepts. Yet, this behavior
is also conditional on the underlying game remaining stationary
throughout the time horizon that it takes players to adapt to it---and
unfortunately, this stationarity assumption is rarely met in practical
applications. In biological models, for example, the reproductive
fitness of an individual may be affected by the ever-changing weather
conditions; in networks, communication channels carry time-dependent
noise and interference as well as signals; and when players try to
sample their strategies, they might have to deal with erroneous or
imprecise readings.

It is thus logical to ask: \textit{does rational behavior still emerge in
the presence of stochastic perturbations that interfere with the
underlying game}?

In evolutionary games, these perturbations traditionally take the form
of ``aggregate shocks'' that are applied directly to the population of
each phenotype. This approach by Fudenberg and Harris \cite{FH92} has
spurred quite a bit of interest and there is a number of features that
differentiate it from the deterministic one. For example, Cabrales
showed in \cite{Ca00} that dominated strategies indeed become extinct,
but only if the variance of the shocks is low enough. More recently,
the work of Imhof and Hofbauer \cite{Im05,Im09} revealed that even
equilibrial play arises over time but again, conditionally on the
variance of the shocks.

Be that as it may, if one looks at games with a finite number of
players, it is hardly relevant to consider shocks of this type because
there are no longer any populations to apply them to. Instead, the
stochastic fluctuations should be reflected directly on the stimuli
that incite players to change their strategies: their payoffs. This
leads to a picture which is very different from the evolutionary one
and is precisely the approach that we will be taking.

\subsection*{Outline of results}

In this paper, we analyze the evolution of players in stochastically
perturbed games of this sort. The particular stimulus-response model
that we consider is simple enough: players keep cumulative scores of
their strategies' performance and employ exponentially more often the
one that scores better. After a few preliminaries in Section \ref
{sec:preliminaries}, this approach is made precise in Section~\ref
{sec:replicator} where we derive the stochastic replicator equation
that governs the behavior of players when their learning curves are
subject to random perturbations.

The replicator equation that we get is different from the ``aggregate
shocks'' approach of \cite{FH92,Ca00,Im05,Im09} and, as a result, it
exhibits markedly different rationality properties as well. In stark
contrast to the results of \cite{Ca00,Im05}, we show in Section \ref
{sec:dominated} that dominated strategies become extinct irrespective
of the noise level (Proposition \ref{prop:dominated}) and provide an
exponential bound for the rate of decay of these strategies
(Proposition \ref{prop:timetolive}). In fact, by induction on the
rounds of elimination of dominated strategies, we show that this is
true even for \textit{iteratively} dominated strategies: despite the
noise, only rationally admissible strategies can survive in the long
run (Theorem \ref{thm:rational}). Then, as an easy corollary of the
above, we infer that players will converge to a strict equilibrium
(Corollary \ref{cor:dominance}) whenever the underlying game is
dominance-solvable.

We continue with the issue of equilibrial play in Section \ref
{sec:congestion} by making a suggestive detour in the land of
congestion games.
If the noise is relatively mild with respect to the rate with which
players learn, we find that the game's potential is a Lyapunov function
which ensures that strict equilibria are stochastically attracting; and
if the game is dyadic (i.e., players only have two choices), this
tameness assumption can be dropped altogether.

Encouraged by the results of Section \ref{sec:congestion}, we attack
the general case in Section~\ref{sec:equilibrium}. As it turns out,
strict equilibria are \textit{always} asymptotically stochastically
stable in the perturbed replicator dynamics that stem from exponential
learning (Theorem~\ref{thm:stability}). This begs to be compared to the
results of \cite{Im05,Im09} where it is the equilibria of a suitably
modified game that are stable, and not necessarily those of the actual
game being played. Fortunately, exponential learning seems to give
players a clearer picture of the original game and there is no need for
similar modifications in our case.

\subsection*{Notational conventions}

Given a finite set $S = \{s_{0},\ldots, s_{n}\}$, we will routinely
identify the set $\Delta(S)$ of probability measures on $S$ with the
standard $n$-dimensional simplex of ${\mathbb R}^{n+1}\dvtx\Delta(S)
\equiv
\{x\in{\mathbb R}^{n+1}\dvtx\sum_{\alpha} x_{\alpha} = 1$ and
$x_{\alpha}\geq0\}$. Under this identification, we will also make no
distinction between $s_{\alpha}\in S$ and the vertex $e_{\alpha}$ of
$\Delta(S)$; in fact, to avoid an overcluttering of indices, we will
frequently use $\alpha$ to refer to either $s_{\alpha}$ or
$e_{\alpha}$, writing, for example, ``$\alpha\in S$'' or ``$u(\alpha)$'' instead
of ``$s_{\alpha}\in S$'' or ``$u(e_{\alpha})$,'' respectively.

To streamline our presentation, we will consistently employ Latin
indices for players ($i,j,k,\ldots$) and Greek for their strategies
($\alpha,\beta,\mu,\ldots$), separating the two by a comma when it would
have been \ae sthetically unpleasant not to. In like manner, when we
have to discriminate between strategies, we will assume that indices
from the first half of the Greek alphabet start at $0$ ($\alpha,\beta=
0,1,2,\ldots$) while those taken from the second half start at $1$
($\mu,\nu= 1,2,\ldots$).

Finally, if $X(t)$ is some stochastic process in ${\mathbb R}^{n}$
starting at
$X(0) = x$, its law will be denoted by $P_{X;x}$ or simply by $P_{x}$
if there is no danger of confusion; and if the context leaves no doubt
as to which process we are referring to, we will employ the term
``almost surely'' in place of the somewhat unwieldy ``$P_{x}$-almost
surely.''

\eject
\section{Preliminaries}
\label{sec:preliminaries}

\subsection{Basic facts and definitions from game theory}

As is customary, our starting point will be a (finite) set of $N$
\textit
{players}, indexed by $i\in\mathcal{N}= \{1,\ldots, N\}$. The players'
possible actions are drawn from their \textit{strategy sets} $\mathcal{S}_{i}
= \{s_{i\alpha}\dvtx\alpha= 0,\ldots, S_{i}-1\}$ and they can
combine them
by choosing their $\alpha_{i}$th (pure) strategy with probability
$p_{i\alpha_{i}}$. In that case, the players' \textit{mixed strategies}
will be described by the points $p_{i} = (p_{i,0},p_{i,1},\ldots)\in
\Delta_{i}:=\Delta(\mathcal{S}_{i})$ or, more succinctly, by the
\textit
{strategy profile} $p = (p_{1},\ldots, p_{N}) \in\Delta:= \prod
_{i}\Delta_{i}$.

In particular, if $e_{i\alpha}$ denotes the $\alpha$th vertex of the
$i$th component simplex \mbox{$\Delta_{i}\hookrightarrow\Delta$}, the (pure)
profile $q = (e_{1,\alpha_{1}},\ldots, e_{N,\alpha_{N}})$ simply
corresponds to player $i$ playing $\alpha_{i}\in\mathcal{S}_{i}$. On the
other hand, if we wish to focus on the strategy of a particular player
$i\in\mathcal{N}$ against that of his \textit{opponents} $\mathcal{N}
_{-i}:=\mathcal{N}
\setminus\{i\}$, we will employ the shorthand notation $(p_{-i};q_{i})
= (p_{1}\cdots q_{i} \cdots p_{N})$ to denote the profile where $i$
plays $q_{i}\in\Delta_{i}$ against his opponents' strategy
$p_{-i}\in
\Delta_{-i}:= \prod_{j\neq i} \Delta_{j}$.

So, once players have made their strategic choices, let $u_{i,\alpha
_{1},\ldots,\alpha_{N}}$ be the reward of player $i$ in the profile
$(\alpha_{1},\ldots,\alpha_{N})\in\mathcal{S}= \prod_{i}\mathcal
{S}_{i}$, that is,
the payoff that strategy $\alpha_{i}\in\mathcal{S}_{i}$ yields to
player $i$
against the strategy $\alpha_{-i}\in\mathcal{S}_{-i}=\prod_{j\neq
i}\mathcal{S}
_{j}$ of $i$'s opponents. Then, if players mix their strategies, their
expected reward will be given by the (multilinear) \textit{payoff
functions} $u_{i}\dvtx\Delta\to{\mathbb R}$:
%
\begin{equation}
u_{i}(p) = \sum_{\alpha_{1}\in\mathcal{S}_{1}} \cdots\sum_{\alpha
_{N}\in
\mathcal{S}_{N}}
u_{i,\alpha_{1}\cdots\alpha_{N}} p_{1,\alpha_{1}} \cdots
p_{N,\alpha_{N}}.
\end{equation}
Under this light, the payoff that a player receives when playing a pure
strategy $\alpha\in\mathcal{S}_{i}$ deserves special mention and
will be
denoted by
%
\begin{equation}
u_{i\alpha}(p) := u_{i}(p_{-i};\alpha) \equiv u_{i}(p_{1}\cdots
\alpha
\cdots p_{N}).
\end{equation}

This collection of \textit{players} $i\in\mathcal{N}$, their \textit
{strategies}
$\alpha_{i}\in\mathcal{S}_{i}$ and their \textit{payoffs} $u_{i}$ will
be our
working definition for a \textit{game in normal form}, usually denoted by
$\mathfrak{G}$---or $\mathfrak{G}(\mathcal{N},\mathcal{S},u)$ if
we need to keep track of
more data.

Needless to say, rational players who seek to maximize their individual
payoffs will avoid strategies that always lead to diminished payoffs
against any play of their opponents. We will thus say that the strategy
$q_{i}\in\Delta_{i}$ is (\textit{strictly}) \textit{dominated} by
$q_{i}'\in
\Delta_{i}$ and we will write $q_{i}\prec q_{i}'$ when
%
\begin{equation}
\label{eq:dominated}
u_{i}(p_{-i};q_{i})<u_{i}(p_{-i};q_{i}')
\end{equation}
for all strategies $p_{-i}\in\Delta_{-i}$ of $i$'s opponents
$\mathcal{N}_{-i}$.

With this in mind, dominated strategies can be effectively removed from
the analysis of a game because rational players will have no incentive
to ever use them. However, by deleting such a strategy, another
strategy (perhaps of another player) might become dominated and further
deletions of \textit{iteratively dominated} strategies might be in order
(see Section \ref{sec:dominated} for more details). Proceeding ad
infinitum, we will say that a strategy is \textit{rationally admissible}
if it survives every round of elimination of dominated strategies. If
the set of rationally admissible strategies is a singleton (e.g., as in
the Prisoner's Dilemma), the game will be called \textit
{dominance-solvable} and the sole surviving strategy will be the game's
\textit{rational solution}.

Then again, not all games can be solved in this way and it is natural
to look for strategies which are stable at least under unilateral
deviations. Hence, we will say that a strategy profile $p\in\Delta$
is a \textit{Nash equilibrium} of the game $\mathfrak{G}$ when
%
\begin{equation}
\label{eq:Nash}
u_{i}(p) \geq u_{i}(p_{-i};q) \qquad\mbox{for all } q\in\Delta_{i},
i\in \mathcal{N}.
\end{equation}
If the equilibrium profile $p$ only contains pure strategies $\alpha
_{i}\in\mathcal{S}_{i}$, we will refer to it as a \textit{pure equilibrium};
and if the inequality (\ref{eq:Nash}) is strict for all $q\neq
p_{i}\in
\Delta_{i},i\in\mathcal{N}$, the equilibrium $p$ will carry instead the
characterization \textit{strict}.

Clearly, if two pure strategies $\alpha,\beta\in\mathcal{S}_{i}$
are present
with positive probability in an equilibrial strategy $p_{i}\in\Delta
_{i}$, then we must have $u_{i\alpha}(p) = u_{i\beta}(p)$ as a result
of $u_{i}$ being linear in $p_{i}$. Consequently, only pure profiles
can satisfy the strict version of (\ref{eq:Nash}) so that strict
equilibria must also be pure. The converse implication is false but
only barely so: a pure equilibrium fails to be strict only if a player
has more than one pure strategies that return the same rewards. Since
this is almost always true (in the sense that the degenerate case can
be resolved by an arbitrarily small perturbation of the payoff
functions), we will relax our terminology somewhat and use the two
terms interchangeably.

To recover the connection of equilibrial play with strategic dominance,
note that if a game is solvable by iterated elimination of dominated
strategies, the single rationally admissible strategy that survives
will be the game's unique strict equilibrium. But the significance of
strict equilibria is not exhausted here: strict equilibria are exactly
the evolutionarily stable strategies of multi-population evolutionary
games---Proposition 5.1 in \cite{We95}. Moreover, as we shall see a
bit later, they are the only asymptotically stable states of the
multi-population replicator dynamics---again, see Chapter 5, pages
216 and 217 of \cite{We95}.

Unfortunately, strict equilibria do not always exist,
Rock-Paper-Scissors being the typical counterexample. Nevertheless,
pure equilibria do exist in many large and interesting classes of
games, even when we leave out dominance-solvable ones. Perhaps the most
noteworthy such class is that of \textit{congestion games}.
\begin{definition}
\label{def:congestion}
A game $\mathfrak{G}\equiv\mathfrak{G}(\mathcal{N},\mathcal{S},u)$
will be called a \textit
{congestion game} when:
\begin{enumerate}
\item all players $i\in\mathcal{N}$ share a common set of \textit{facilities}
$\mathcal{F}$ as their strategy set: $\mathcal{S}_{i}=\mathcal{F}$
for all
$i\in\mathcal{N}$;
\item the payoffs are functions of the number of players sharing a
particular facility: $u_{i,\alpha_{1}\cdots\alpha\cdots\alpha_{N}}
\equiv u_{\alpha}(N_{\alpha})$ where $N_{\alpha}$ is the number of
players choosing the same facility as $i$.
\end{enumerate}
\end{definition}

Amazingly enough, Monderer and Shapley made the remarkable discovery in
\cite{MS96} that these games are actually equivalent to the class of
\textit{potential games}.
\begin{definition}
\label{def:potential}
A game $\mathfrak{G}\equiv\mathfrak{G}(\mathcal{N},\mathcal{S},u)$
will be called a \textit
{potential game} if there exists a function $V\dvtx\Delta\to{\mathbb R}$
such that
%
\begin{equation}
\label{eq:potential}
u_{i}(p_{-i};q_{i}) - u_{i}(p_{-i};q'_{i}) = - \bigl(V(p_{-i};q_{i}) -
V(p_{-i};q_{i}') \bigr)
\end{equation}
for all players $i\in\mathcal{N}$ and all strategies $p_{-i}\in
\Delta
_{-i}$, $q_{i},q_{i}'\in\Delta_{i}$.
\end{definition}

This equivalence reveals that both classes of games possess equilibria
in pure strategies: it suffices to look at the vertices of the face of
$\Delta$ where the (necessarily multilinear) potential function $V$
is minimized.

\subsection{Learning, evolution and the replicator dynamics}
As one would expect, locating the Nash equilibria of a game is a rather
complicated problem that requires a great deal of global calculations,
even in the case of potential games (where it reduces to minimizing a
multilinear function over a convex polytope). Consequently, it is of
interest to see whether there are simple and distributed learning
schemes that allow players to arrive at a reasonably stable solution.

One such scheme is based on an exponential learning behavior where
players play the game repeatedly and keep records of their strategies'
performance. In more detail, at each instance of the game all players
$i\in\mathcal{N}$ update the cumulative scores $U_{i\alpha}$ of their
strategies $\alpha\in\mathcal{S}_{i}$ as specified by the recursive formula
%
\begin{equation}
\label{eq:discretescore}
U_{i\alpha} (t+1) = U_{i\alpha}(t) + u_{i\alpha}(p(t)),
\end{equation}
where $p(t)\in\Delta$ is the players' strategy profile at the $t$th
iteration of the game and, in the absence of initial bias, we assume
that $U_{i\alpha}(0) = 0$ for all $i\in\mathcal{N}, \alpha\in
\mathcal{S}_{i}$.
These scores reinforce the perceived success of each strategy as
measured by the average payoff it yields and hence, it stands to reason
that players will lean towards the strategy with the highest score. The
precise way in which they do that is by playing according to the
namesake exponential law:
%
\begin{equation}
\label{eq:discretelogit}
p_{i\alpha}(t+1) = \frac{e^{U_{i\alpha}(t+1)}}{\sum_{\beta\in
\mathcal{S}
_{i}}e^{U_{i\beta}(t+1)}}.
\end{equation}

For simplicity, we will only consider the case where players update
their scores in continuous time, that is, according to the coupled equations
%
\begin{subequation}
\begin{eqnarray}
\label{eq:detscore}
dU_{i\alpha}(t) &=& u_{i\alpha} (x(t)) \,dt,\\
\label{eq:detlogit}
x_{i\alpha}(t) &=& \frac{e^{U_{i\alpha}(t)}}{\sum_{\beta
}e^{U_{i\beta}(t)}}.
\end{eqnarray}
\end{subequation}
Then, if we differentiate (\ref{eq:detlogit}) to decouple it from
(\ref
{eq:detscore}), we obtain the \textit{standard} (\textit{multi-population})
\textit{replicator dynamics}
%
\begin{equation}
\label{eq:RD}
\frac{dx_{i\alpha}}{dt} =
x_{i\alpha} \biggl(u_{i\alpha}(x) - \sum_{\beta} x_{i\beta} u_{i\beta
}(x) \biggr) =
x_{i\alpha} \bigl(u_{i\alpha}(x) - u_{i}(x) \bigr).
\end{equation}

Alternatively, if players learn at different speeds as a result of
varied stimulus-response characteristics, their updating will take the form
%
\begin{equation}
\label{eq:ratelogit}
x_{i\alpha}(t) = \frac{e^{\lambda_{i} U_{i\alpha}(t)}}{\sum_{\beta}
e^{\lambda_{i} U_{i\beta}(t)}},
\end{equation}
where $\lambda_{i}$ represents the \textit{learning rate} of player $i$,
that is, the ``weight'' which he assigns to his perceived scores
$U_{i\alpha}$. In this way, the replicator equation evolves at a
different time scale for each player, leading to the \textit
{rate-adjusted} dynamics
%
\begin{equation}
\label{eq:LRD}
\frac{dx_{i\alpha}}{dt} = \lambda_{i} x_{i\alpha} \bigl(u_{i\alpha}(x) -
u_{i}(x) \bigr).
\end{equation}
Naturally, the uniform dynamics (\ref{eq:RD}) are recovered when all
players learn at the ``standard'' rate $\lambda_{i}=1$.

If we view the exponential learning model (\ref{eq:discretelogit}) from
a stimulus-response angle, we see that that the payoff of a strategy
simply represents an (exponential) propensity of employing said
strategy. It is thus closely related to the algorithm of \textit{logistic
fictitious play} \cite{FL98} where the strategy $x_{i}$ of (\ref
{eq:ratelogit}) can be seen as the (unique) best reply to the profile
$x_{-i}$ in some suitably modified payoffs\vspace*{1pt} $v_{i}(x) = u_{i}(x) + \frac
{1}{\lambda_{i}}H(x_{i})$. Interestingly enough, $H(x_{i})$ turns out
to be none other than the \textit{entropy} of $x_{i}$:
%
\begin{equation}
H(x_{i}) = - \sum_{\beta\dvtx x_{i\beta}>0} x_{i\beta} \log
x_{i\beta}.
\end{equation}
That being so, we deduce that the learning rates $\lambda_{i}$ act the
part of (player-specific) inverse temperatures: in high temperatures
(small $\lambda_{i}$), the players' learning curves are ``soft'' and
the payoff differences between strategies are toned down; on the
contrary, if $\lambda_{i}\to\infty$ the scheme ``freezes'' to a myopic
best-reply process.

The replicator dynamics were first derived in \cite{TJ78} in the
context of population biology, first for different phenotypes within a
single species (single-population models), and then for different
species altogether (multi-population models; \cite
{HS88} and \cite{We95} provide excellent surveys). In both these cases, one begins with
large populations of individuals that are programmed to a particular
behavior (e.g., \textsf{fight} for ``hawks'' or \textsf{flight} for
``doves'') and matches them randomly in a game whose payoffs directly
affect the reproductive fitness of the individual players.

More precisely, let $z_{i\alpha}(t)$ be the population size of the
phenotype (strategy) $\alpha\in\mathcal{S}_{i}$ of species (player)
$i\in
\mathcal{N}$ in some multi-population model where individuals are
matched to
play a game $\mathfrak{G}$ with payoff functions $u_{i}$. Then, the relative
frequency (share) of $\alpha$ will be specified by the \textit{population
state} $x=(x_{1},\ldots, x_{N})\in\Delta$ where $x_{i\alpha} =
z_{i\alpha}/\sum_{\beta}z_{i\beta}$. So, if $N$ individuals are drawn
randomly from the $N$ species, their expected payoffs will be given by
$u_{i}(x)$, $i\in\mathcal{N}$, and if these payoffs represent a proportionate
increase in the phenotype's fitness (measured as the number of
offsprings in the unit of time), we will have
%
\begin{equation}
\label{eq:offspring}
d z_{i\alpha}(t) = z_{i\alpha}(t) u_{i\alpha}(x(t)) \,dt.
\end{equation}
As a result, the population state $x(t)$ will evolve according to
%
\begin{equation}
\label{eq:ERD}\quad
\frac{d x_{i\alpha}}{dt} = \frac{1}{\sum_{\beta} z_{i\beta}}
\,\frac{d
z_{i\alpha}}{d t}
- \sum_{\gamma} \frac{x_{i\alpha}}{\sum_{\beta} z_{i\beta}}
\,\frac{d
z_{i\gamma}}{dt}
= x_{i\alpha} \bigl(u_{i\alpha}(x) - u_{i}(x) \bigr),
\end{equation}
which is exactly (\ref{eq:RD}) viewed from an evolutionary perspective.

On the other hand, we should note here that in \textit{single-population}
models the resulting equation is cubic and not quadratic because
strategies are matched against themselves. To wit, assume that
individuals are randomly drawn from a large population and are matched
against one another in a (symmetric) 2-player game $\mathfrak{G}$ with
strategy space $\mathcal{S}=\{1,\ldots, S\}$ and payoff matrix $u=\{
u_{\alpha
\beta}\}$. Then, if $x_{\alpha}$ denotes the population share of
individuals that are programmed to the strategy $\alpha\in\mathcal{S}$,
their expected payoff in a random match will be given by $u_{\alpha}(x)
:= \sum_{\beta} u_{\alpha\beta}x_{\beta} \equiv u(\alpha,x)$;
similarly, the population average payoff will be $u(x,x) = \sum
_{\alpha
} x_{\alpha} u_{\alpha}(x)$. Hence, by following the same procedure as
above, we end up with the single-population replicator dynamics
%
\begin{equation}
\label{eq:1RD}
\frac{d x_{\alpha}}{dt} = x_{\alpha} \bigl(u_{\alpha}(x) - u(x,x)
\bigr),
\end{equation}
which behave quite differently than their multi-population counterpart
(\ref{eq:ERD}).

As far as rational behavior is concerned, the replicator dynamics have
some far-reaching ramifications. If we focus on multi-population
models, Samuelson and Zhang showed in \cite{SZ92} that the share
$x_{i\alpha}(t)$ of a strategy $\alpha\in\mathcal{S}_{i}$ which is strictly
dominated (even iteratively) converges to zero along any interior
solution path of~(\ref{eq:RD}); in other words, \textit{dominated
strategies become extinct in the long run}. Additionally, there is a
remarkable equivalence between the game's Nash equilibria and the
stationary points of the replicator dynamics: \textit{the asymptotically
stable states of \textup{(\ref{eq:RD})} coincide precisely with the strict Nash
equilibria of the underlying game}~\cite{We95}.

\subsection{Elements of stability analysis}

A large part of our work will be focused on examining whether the
rationality properties of exponential learning (elimination of
dominated strategies and asymptotic stability of strict equilibria)
remain true in a stochastic setting. However, since asymptotic
stability is (usually) too stringent an expectation for stochastic
dynamical systems, we must instead consider its stochastic analogue.

That being the case, let $W(t)=(W_{1}(t),\ldots, W_{n}(t))$ be a standard
Wiener process in ${\mathbb R}^{n}$ and consider the stochastic differential
equation (SDE)
%
\begin{equation}
\label{eq:SDE}
dX_{\alpha}(t) = b_{\alpha}(X(t))\,dt + \sum_{\beta} \sigma
_{\alpha
\beta}(X(t))\, dW_{\beta}(t).
\end{equation}
Following \cite{GS71,Ar74}, the notion of asymptotic stability in this
SDE is expressed by the following.
\begin{definition}
\label{def:stability}
We will say that $q\in{\mathbb R}^{n}$ is \textit{stochastically
asymptotically
stable} when, for every neighborhood $U$ of $q$ and every $\varepsilon>0$,
there exists a neighborhood $V$ of $q$ such that
%
\begin{equation}
P_{x} \Bigl\{X(t)\in U\mbox{ for all }t\geq0, \lim_{t\to\infty} X(t) =
q \Bigr\}
\geq1-\varepsilon
\end{equation}
for all initial conditions $X(0) = x\in V$ of the SDE (\ref{eq:SDE}).
\end{definition}

Much the same as in the deterministic case, stochastic asymptotic
stability is often established by means of a Lyapunov function. In our
context, this notion hinges on the second order differential operator
that is associated to (\ref{eq:SDE}), namely the \textit
{generator} $L$ of $X(t)$:
%
\begin{equation}
L= \sum_{\alpha=1}^{n} b_{\alpha}(x) \,\frac{\partial}{\partial
x_{\alpha}}
+ \frac{1}{2} \sum_{\alpha,\beta=1}^{n} (\sigma(x)\sigma^{T}(x)
)_{\alpha\beta}\,
\frac{\partial^{2}}{\partial x_{\alpha}\, \partial x_{\beta}}.
\end{equation}
The importance of this operator can be easily surmised from It\^{o}'s
lemma; indeed, if $f\dvtx{\mathbb R}^{n}\to{\mathbb R}$ is
sufficiently smooth, the generator
$L$ simply captures the drift of the process $Y(t) = f(X(t))$:
%
\begin{equation}
dY(t) = Lf(X(t))\,dt + \sum_{\alpha,\beta} \frac{\partial f}{\partial
x_{\alpha}} \bigg|_{X(t)} \sigma_{\alpha\beta}(X(t))\,dW_{\beta}(t).
\end{equation}
In this way, $L$ can be seen as the stochastic version of the time
derivative $\frac{d}{dt}$; this analogy then leads to the following.
\begin{definition}
\label{def:Lyapunov}
Let $q\in{\mathbb R}^{n}$ and let $U$ be an open neighborhood of $q$.
We will
say that $f$ is a (local) \textit{stochastic Lyapunov function} for the
SDE (\ref{eq:SDE}) if:
\begin{enumerate}
\item$f(x)\geq0$ for all $x\in U$, with equality iff $x=q$;
\item there exists a constant $k>0$ such that $Lf(x) \leq- k f(x)$
for all $x\in U$.
\end{enumerate}
\end{definition}

Whenever such a Lyapunov function exists, it is known that the point
$q\in{\mathbb R}^{n}$ where $f$ attains its minimum will be stochastically
asymptotically stable---for example, see Theorem 4 in pages 314 and 315 of
\cite{GS71}.

A final point that should be mentioned here is that our analysis will
be constrained on the compact polytope $\Delta=\prod_{i}\Delta_{i}$
instead of all of $\prod_{i}{\mathbb R}^{S_{i}}$. Accordingly, the
``neighborhoods'' of Definitions \ref{def:stability} and \ref
{def:Lyapunov} should be taken to mean ``neighborhoods in $\Delta$,''
that is, neighborhoods in the subspace topology of $\Delta
\hookrightarrow
\prod_{i}{\mathbb R}^{S_{i}}$. This minor point should always be clear
from the
context and will only be raised in cases of ambiguity.

\section{Learning in the presence of noise}
\label{sec:replicator}

Of course, it could be argued that the rationality properties of the
exponential learning scheme are a direct consequence of the players'
receiving accurate information about the game when they update their
scores. However, this is a requirement that cannot always be met: the
interference of nature in the game or imperfect readings of one's
utility invariably introduce fluctuations in (\ref{eq:detscore}), and
in their turn, these lead to a perturbed version of the replicator
dynamics (\ref{eq:RD}).

To account for these random perturbations, we will assume that the
players' scores are now governed instead by the \textit{stochastic}
differential equation
%
\begin{equation}
\label{eq:score}
dU_{i\alpha} (t) = u_{i\alpha}(X(t)) \,dt + \eta_{i\alpha}(X(t))\,
dW_{i\alpha}(t),
\end{equation}
where, as before, the strategy profile $X(t)\in\Delta$ is given by
the logistic law
%
\begin{equation}
\label{eq:logit}
X_{i\alpha}(t) = \frac{e^{U_{i\alpha}(t)}}{\sum_{\beta}
e^{U_{i\beta}(t)}}.
\end{equation}
In this last equation, $W(t)$ is a standard Wiener process living in
$\prod_{i}{\mathbb R}^{S_{i}}$ and the coefficients $\eta_{i\alpha}$
measure the
impact of the noise on the players' scoring systems. Of course, these
coefficients need not be constant: after all, the effect of the noise
on the payoffs might depend on the state of the game in some typically
continuous way. For this reason, we will assume that the functions
$\eta
_{i\alpha}$ are continuous on $\Delta$, and we will only note en
passant that our results still hold for essentially bounded
coefficients $\eta_{i\alpha}$ (we will only need to replace $\min$ and
$\max$ with $\operatorname{ess}\inf$ and $\operatorname{ess}\sup$,
respectively, in all expressions
involving $\eta_{i\alpha}$).

A very important instance of this dependence can be seen if
$\eta_{i\alpha}(x_{-i};\alpha)=0$ for all $i\in\mathcal{N},
\alpha\in\mathcal{S}_{i}, x_{-i}\in\Delta_{-i}$, in which case equation
(\ref{eq:score}) becomes a convincing model for the case of
insufficient information. It states that when a player actually uses a
strategy, his payoff observations are accurate enough; but with regards
to strategies he rarely employs, his readings could be arbitrarily off
the mark.

Now, to decouple (\ref{eq:score}) and (\ref{eq:logit}), we may simply
apply It\^{o}'s lemma to the process $X(t)$. To that end, recall that
$W(t)$ has independent components across players and strategies, so
that $dW_{j\beta}\cdot dW_{k\gamma} = \delta_{jk}\delta_{\beta
\gamma
}\,dt$ (the Kronecker symbols $\delta_{\beta\gamma}$ being $0$ for
$\beta\neq\gamma$ and $1$, otherwise). Then, It\^{o}'s formula gives
%
\begin{eqnarray}
\label{eq:dXprelim}
dX_{i\alpha} &=& \sum_{j}\sum_{\beta} \frac{\partial X_{i\alpha
}}{\partial
U_{j\beta}} \,dU_{j\beta}\nonumber\\
&&{} + \frac{1}{2}\sum_{j,k}\sum_{\beta
,\gamma
} \frac{\partial^{2} X_{i\alpha}}{\partial U_{j\beta}\,\partial
U_{k\gamma}} \,d
U_{j\beta}\cdot dU_{k\gamma}\nonumber\\[-8pt]\\[-8pt]
&=& \sum_{\beta} \biggl(u_{i\beta}(X)\,\frac{\partial X_{i\alpha}}{\partial
U_{i\beta}}
+ \frac{1}{2}\eta_{i\beta}^{2}(X)\,\frac{\partial^{2}X_{i\alpha
}}{\partial
U_{i\beta}^{2}} \biggr) \,dt\nonumber\\
&&{} + \sum_{\beta} \eta_{i\beta}(X)\frac{\partial X_{i\alpha
}}{\partial
U_{i\beta
}} \,dW_{i\beta}.\nonumber
\end{eqnarray}
On the other hand, a simple differentiation of (\ref{eq:logit}) yields
%
\begin{subequation}
\begin{eqnarray}
\frac{\partial X_{i\alpha}}{\partial U_{i\beta}} &=& X_{i\alpha
}(\delta
_{\alpha
\beta} - X_{i\beta}),\\
\frac{\partial^{2} X_{i\alpha}}{\partial U_{i\beta}^{2}} &=&
X_{i\alpha
}(\delta
_{\alpha\beta} - X_{i\beta})(1-2X_{i\beta})
\end{eqnarray}
\end{subequation}
and by plugging these expressions back into (\ref{eq:dXprelim}), we get
%
\begin{eqnarray}
\label{eq:SRD}\qquad
dX_{i\alpha} &=& X_{i\alpha}
[u_{i\alpha}(X) - u_{i}(X) ] \,dt\nonumber\\
&&{} + X_{i\alpha} \biggl[\frac{1}{2}\eta_{i\alpha}^{2}(X)(1-2 X_{i\alpha}) -
\frac{1}{2}\sum_{\beta}\eta_{i\beta}^{2}(X) X_{i\beta
}(1-2X_{i\beta})
\biggr] \,dt\\
&&{} + X_{i\alpha} \biggl[\eta_{i\alpha}(X) \,dW_{i\alpha} - \sum_{\beta
} \eta
_{i\beta}(X) X_{i\beta} \,dW_{i\beta} \biggr]\nonumber.
\end{eqnarray}
Alternatively, if players update their strategies with different
learning rates $\lambda_{i}$, we should instead apply It\^{o}'s formula
to (\ref{eq:ratelogit}). In so doing, we obtain
\renewcommand{\theequation}{3.5$'$}
\begin{eqnarray}
\label{eq:SLRD}\qquad
dX_{i\alpha} &=& \lambda_{i}X_{i\alpha}
[u_{i\alpha}(X) - u_{i}(X) ] \,dt\nonumber\\
&&{} +\frac{\lambda_{i}^{2}}{2}X_{i\alpha} \biggl[\eta_{i\alpha}^{2}(X)(1-2
X_{i\alpha}) -\sum_{\beta}\eta_{i\beta}^{2}(X) X_{i\beta
}(1-2X_{i\beta
}) \biggr] \,dt\nonumber\\[-8pt]\\[-8pt]
&&{} +\lambda_{i} X_{i\alpha} \Bigl[\eta_{i\alpha}(X)\,
d
W_{i\alpha} - \sum\eta_{i\beta}(X) X_{i\beta} \,dW_{i\beta}
\Bigr]\nonumber\\
&=& b_{i\alpha}(X) \,dt + \sum_{\beta} \sigma_{i,\alpha\beta
}(X) \,d
W_{i\beta},\nonumber
\end{eqnarray}
where, in obvious notation, $b_{i\alpha}(x)$ and $\sigma_{i,\alpha
\beta
}(x)$ are, respectively, the drift and diffusion coefficients of the
diffusion $X(t)$. Obviously, when $\lambda_{i}=1$, we recover the
uniform dynamics (\ref{eq:SRD}); equivalently (and this is an
interpretation that is well worth keeping in mind), the rates $\lambda
_{i}$ can simply be regarded as a commensurate inflation of the payoffs
and noise coefficients of player $i\in\mathcal{N}$ in the uniform logistic
model (\ref{eq:logit}).

Equation (\ref{eq:SRD}) and its rate-adjusted sibling (\ref{eq:SLRD})
will constitute our stochastic version of the replicator dynamics and
thus merit some discussion in and by themselves. First, note that these
dynamics admit a (unique) strong solution for any initial state $X(0) =
x\in\Delta$, even though they do not satisfy the linear growth
condition $|b(x)|+|\sigma(x)|\leq C(1+|x|)$ that is required for the
existence and uniqueness theorem for SDEs (e.g., Theorem 5.2.1 in \cite
{Ok07}). Instead, an addition over $\alpha\in\mathcal{S}_{i}$
reveals that
every simplex $\Delta_{i}\subseteq\Delta$ remains invariant under
(\ref{eq:SRD}): if $X_{i}(0) = x_{i}\in\Delta_{i}$, then $d
(\sum
_{\alpha}X_{i\alpha} ) = 0$ and hence, $X_{i}(t)$ will stay in
$\Delta
_{i}$ for all $t\geq0$---actually, it is not harder to see that every
face of $\Delta$ is a trap for $X(t)$.

So, if $\phi$ is a smooth bump function that is equal to $1$ on some
open neighborhood of $U\supseteq\Delta$ and which vanishes outside
some compact set $K\supseteq U$, the SDE
%
\setcounter{equation}{5}
\renewcommand{\theequation}{\arabic{section}.\arabic{equation}}
\begin{equation}
dX_{i\alpha} = \phi(X) \biggl(b_{i\alpha}(X) \,dt + \sum_{\beta}
\sigma
_{i,\alpha\beta}(X) \,dW_{i\beta} \biggr)
\end{equation}
\textit{will} have bounded diffusion and drift coefficients and will thus
admit a unique strong solution. But since this last equation agrees
with (\ref{eq:SRD}) on $\Delta$ and any solution of (\ref{eq:SRD})
always stays in $\Delta$, we can easily conclude that our perturbed
replicator dynamics admit a unique strong solution for any initial
$X(0) = x\in\Delta$.

It is also important to compare the dynamics (\ref{eq:SRD}), (\ref
{eq:SLRD}) to the ``aggregate shocks'' approach of Fudenberg and Harris
\cite{FH92} that has become the principal incarnation of the replicator
dynamics in a stochastic environment. So, let us first recall how
aggregate shocks enter the replicator dynamics in the first place. The
main idea is that the reproductive fitness of an individual is not only
affected by deterministic factors but is also subject to stochastic
shocks due to the ``weather'' and the interference of nature with the
game. More precisely, if $Z_{i\alpha}(t)$ denotes the population size
of phenotype $\alpha\in\mathcal{S}_{i}$ of the species $i\in
\mathcal{N}$ in some
multi-population evolutionary game $\mathfrak{G}$, its growth will be
determined by
%
\begin{equation}
dZ_{i\alpha}(t) = Z_{i\alpha}(t) \bigl(u_{i\alpha}(X(t))\,dt + \eta
_{i\alpha}\,dW_{i\alpha}(t) \bigr),
\end{equation}
where, as in (\ref{eq:offspring}), $X(t)\in\Delta$ denotes the
population shares $X_{i\alpha} = Z_{i\alpha}/\sum_{\beta}Z_{i\beta}$.
In this way, It\^{o}'s lemma yields the \textit{replicator dynamics with
aggregate shocks}:
%
\begin{eqnarray}
\label{eq:ASRD}
dX_{i\alpha} &=& X_{i\alpha}
\biggl[ \bigl(u_{i\alpha}(X) - u_{i}(X) \bigr)
- \biggl(\eta_{i\alpha}^{2} X_{i\alpha} - \sum_{\beta} \eta_{i\beta}^{2}
X_{i\beta}^{2} \biggr) \biggr] \,dt\nonumber\\[-8pt]\\[-8pt]
&&{} + X_{i\alpha} \Bigl[\eta_{i\alpha} \,dW_{i\alpha} - \sum\eta
_{i\beta}
X_{i\beta} \,dW_{i\beta} \Bigr].\nonumber
\end{eqnarray}

We thus see that the effects of noise propagate differently in the case
of exponential learning and in the case of evolution. Indeed, if we
compare equations (\ref{eq:SRD}) and (\ref{eq:ASRD}) term by term, we
see that the drifts are not quite the same: even though the payoff
adjustment $u_{i\alpha} - u_{i}$ ties both equations back together in
the deterministic setting ($\eta=0$), the two expressions differ by
%
\begin{equation}
\label{eq:expterm}
X_{i\alpha} \biggl[\frac{1}{2}\eta_{i\alpha}^{2} - \frac{1}{2}\sum
_{\beta}
\eta_{i\beta}^{2} X_{i\beta} \biggr] \,dt.
\end{equation}
Innocuous as this term might seem, it is actually crucial for the
rationality properties of exponential learning in games with randomly
perturbed payoffs. As we shall see in the next sections, it leads to
some miraculous cancellations that allow rationality to emerge in all
noise levels.

This difference further suggests that we can pass from (\ref{eq:SRD})
to (\ref{eq:ASRD}) simply by modifying the game's payoffs to
$\widetilde
{u}_{i\alpha} = u_{i\alpha} + \frac{1}{2}\eta_{i\alpha}^{2}$. Of
course, this presumes that the noise coefficients $\eta_{i\alpha}$ be
constant---the general case would require us to allow for games whose
payoffs may not be multilinear. This apparent lack of generality does
not really change things but we prefer to keep things simple and for
the time being, it suffices to point out that this modified game was
precisely the one that came up in the analysis of \cite{Im05,Im09}. As
a result, this modification appears to play a pivotal role in setting
apart learning and evolution in a stochastic setting: whereas the
modified game is deeply ingrained in the process of natural selection,
exponential learning seems to give players a clearer picture of the
actual underlying game.

\section{Extinction of dominated strategies}
\label{sec:dominated}

Thereby armed with the stochastic replicator equations (\ref
{eq:SRD}), (\ref{eq:SLRD}) to model exponential learning in noisy
environments, the logical next step is to see if the rationality
properties of the deterministic dynamics carry over to this stochastic
setting. In this direction, we will first show that dominated
strategies always become extinct in the long run and that only the
rationally admissible ones survive.

As in \cite{Ca00} (implicitly) and \cite{Im05} (explicitly), the key
ingredient of our approach will be the \textit{cross entropy} between two
mixed strategies $q_{i}, x_{i}\in\Delta_{i}$ of player $i\in\mathcal{N}$:
%
\begin{equation}
\label{eq:cross}
H(q_{i},x_{i}) := - \sum_{\alpha\dvtx q_{i\alpha}>0} q_{i\alpha}
\log
(x_{i\alpha}) \equiv H(q_{i}) + d_{\operatorname{KL}}(q_{i},x_{i}),
\end{equation}
where $H(q_{i}) = -\sum_{\alpha}q_{i\alpha} \log q_{i\alpha}$ is the
\textit{entropy} of $q_{i}$ and $d_{\operatorname{KL}}$ is the
intimately related
\textit
{Kullback--Leibler divergence} (or \textit{relative entropy}):
%
\begin{equation}
\label{eq:KL}
d_{\operatorname{KL}}(q_{i},x_{i}) := H(q_{i},x_{i}) - H(q_{i}) = \sum
_{\alpha\dvtx
q_{i\alpha}>0} q_{i\alpha} \log\frac{q_{i\alpha}}{x_{i\alpha}}.
\end{equation}
This divergence function is central in the stability analysis of the
(deterministic) replicator dynamics because it serves as a distance
measure in probability space~\cite{We95}. As it stands however,
$d_{\operatorname{KL}}$
is not a distance function per se: neither is it symmetric, nor does it
satisfy the triangle inequality. Still, it has the very useful property
that $d_{\operatorname{KL}}(q_{i},x_{i}) < \infty$ iff $x_{i}$
employs with positive
probability all pure strategies $\alpha\in\mathcal{S}_{i}$ that are present
in $q_{i}$ [i.e., iff $\operatorname{supp}(q_{i}) \subseteq
\operatorname{supp}(x_{i})$ or iff
$q_{i}$ is absolutely continuous w.r.t. $x_{i}$]. Therefore, if
$d_{\operatorname{KL}}
(q_{i},x_{i}) = \infty$ for all dominated strategies $q_{i}$ of player
$i$, it immediately follows that $x_{i}$ cannot be dominated itself. In
this vein, we have the following.
\begin{proposition}
\label{prop:dominated}
Let $X(t)$ be a solution of the stochastic replicator dynamics (\ref
{eq:SRD}) for some interior initial condition $X(0) = x\in
\operatorname{Int}
(\Delta
)$. Then, if $q_{i}\in\Delta_{i}$ is (strictly) dominated,
%
\begin{equation}
\lim_{t\to\infty} d_{\operatorname{KL}}(q_{i},X_{i}(t)) = \infty
\qquad\mbox{almost surely.}
\end{equation}
In particular, if $q_{i}=\alpha\in\mathcal{S}_{i}$ is pure, we will have
$\lim_{t\to\infty} X_{i\alpha}(t) = 0$ (a.s.): strictly dominated
strategies do not survive in the long run.
\end{proposition}
\begin{pf}
Note first that $X(0) = x \in\operatorname{Int}(\Delta)$ and hence,
$X_{i}(t)$ will
almost surely stay in $\operatorname{Int}(\Delta_{i})$ for all $t\geq
0$; this is a
simple consequence of the uniqueness of strong solutions and the
invariance of the faces of $\Delta_{i}$ under the dynamics (\ref{eq:SRD}).

Let us now consider the cross entropy $G_{q_{i}}(t)$ between $q_{i}$
and $X_{i}(t)$:
%
\begin{equation}
G_{q_{i}}(t) \equiv H(q_{i},X_{i}(t)) = -\sum_{\alpha} q_{i\alpha}
\log X_{i\alpha}(t).
\end{equation}
As a result of $X_{i}(t)$ being an interior path, $G_{q_{i}}(t)$ will
remain finite for all $t\geq0$ (a.s.). So, by applying It\^{o}'s lemma
we get
%
\begin{eqnarray}
dG_{q_{i}} &=& \sum_{\beta} \frac{\partial G_{q_{i}}}{\partial
X_{i\beta}}\,dX_{i\beta}
+ \frac{1}{2} \sum_{\beta,\gamma} \frac{\partial
^{2}G_{q_{i}}}{\partial
X_{i\gamma}\,\partial X_{i\beta}} \,dX_{i\beta}\cdot dX_{i\gamma}
\nonumber\\[-8pt]\\[-8pt]
&=& -\sum_{\beta} \frac{q_{i\beta}}{X_{i\beta}} \,dX_{i\beta}
+ \frac{1}{2} \sum_{\beta} \frac{q_{i\beta}}{X_{i\beta}^{2}} (d
X_{i\beta} )^{2}\nonumber
\end{eqnarray}
and, after substituting $dX_{i\beta}$ from the dynamics (\ref
{eq:SRD}), this last equation becomes
%
\begin{eqnarray}
dG_{q_{i}} &=& \sum_{\beta} q_{i\beta} \biggl[ u_{i}(X) - u_{i\beta}(X)
+\frac{1}{2} \sum_{\gamma} \eta_{i\gamma}^{2}(X) X_{i\gamma
}(1-X_{i\gamma}) \biggr] \,dt\nonumber\\[-8pt]\\[-8pt]
&&{} +\sum_{\beta} q_{i\beta}\sum_{\gamma} (X_{i\gamma} - \delta
_{\beta
\gamma}) \eta_{i\gamma} (X)\,dW_{i\gamma}.\nonumber
\end{eqnarray}

Accordingly, if $q_{i}'\in\Delta_{i}$ is another mixed strategy of
player $i$, we readily obtain
%
\begin{eqnarray}
\label{eq:Gcomp}
d G_{q_{i}} - d G_{q'_{i}} &=& \bigl(u_{i}(X_{-i};q'_{i}) -
u_{i}(X_{-i};q_{i}) \bigr) \,dt\nonumber\\[-8pt]\\[-8pt]
&&{} + \sum_{\beta} (q_{i\beta}' -
q_{i\beta})
\eta_{i\beta}(X) \,dW_{i\beta}\nonumber
\end{eqnarray}
and, after integrating,
%
\begin{eqnarray}
\label{eq:Gint}
G_{q_{i}-q'_{i}}(t) &=& H(q_{i}-q_{i}',x)
+\int_{0}^{t} u_{i}\bigl(X_{-i}(s); q_{i}'-q_{i}\bigr)\,ds \nonumber\\[-8pt]\\[-8pt]
&&{} + \sum_{\beta}(q_{i\beta}'-q_{i\beta})\int_{0}^{t}\eta_{i\beta
}(X(s)) \,dW_{i\beta}(s).\nonumber
\end{eqnarray}
Suppose then that $q_{i}\prec q'_{i}$ and let $v_{i} =\inf\{
u_{i}(x_{-i};q'_{i}-q_{i})\dvtx x_{-i}\in\Delta_{-i}\}$. With
$\Delta
_{-i}$ compact, it easily follows that $v_{i}>0$ and the first term of
(\ref{eq:Gint}) will be bounded from below by $v_{i}t$.

However, since monotonicity fails for It\^{o} integrals, the second
term must be handled with more care. To that end, let $\xi_{i}(s) =
\sum
_{\beta} (q'_{i\beta}-q_{i\beta}) \eta_{i\beta}(X(s))$ and note that
the Cauchy--Schwarz inequality gives
%
\begin{eqnarray}
\xi_{i}^{2}(s) &\leq& S_{i} \sum_{\beta}(q'_{i\beta}-q_{i\beta})^{2}
\eta_{i\beta}^{2}(X(s)) \nonumber\\[-8pt]\\[-8pt]
&\leq& S_{i}\eta_{i}^{2} \sum_{\beta
}(q'_{i\beta
}-q_{i\beta})^{2}\leq 2 S_{i} \eta_{i}^{2},\nonumber
\end{eqnarray}
where $S_{i} = |\mathcal{S}_{i}|$ is the number of pure strategies available
to player $i$ and $\eta_{i} = \max\{|\eta_{i\beta}(x)|\dvtx x\in
\Delta
,\beta\in\mathcal{S}_{i}\}$; recall also that $q_{i},q_{i}'\in
\Delta_{i}$
for the last step. Therefore, if $\psi_{i}(t)=\sum_{\beta}
(q'_{i\beta
}-q_{i\beta})\int_{0}^{t}\eta_{i\beta}(X(s))\,dW_{i\beta}(s)$ denotes
the martingale part of (\ref{eq:Gcomp}) and $\rho_{i}(t)$ is its
quadratic variation, the previous inequality yields
%
\begin{equation}
\label{eq:quvar}
\rho_{i}(t) = [\psi_{i},\psi_{i}](t) = \int_{0}^{t}\xi
_{i}^{2}(s)\,ds
\leq2 S_{i} \eta_{i}^{2} t.
\end{equation}

Now, if $\lim_{t\to\infty}\rho_{i}(t)=\infty$, it follows from the
time-change theorem for martingales (e.g., Theorem 3.4.6 in \cite
{KS98}) that there exists a Wiener process $\widetilde{W}_{i}$ such
that $\psi_{i}(t) = \widetilde{W}_{i}(\rho_{i}(t))$. Hence, by the law
of the iterated logarithm we get
%
\begin{eqnarray}
&&\liminf_{t\to\infty} G_{q_{i}-q_{i}'}(t) \nonumber\\
&&\qquad\geq H(q_{i}-q_{i}',x) +
\liminf_{t\to\infty} \bigl(v_{i}t + \widetilde{W}_{i}(\rho_{i}(t))\bigr)\nonumber\\
&&\qquad\geq H(q_{i}-q_{i}',x) + \liminf_{t\to\infty} \bigl(v_{i}t - \sqrt
{2\rho_{i}(t)\log\log\rho_{i}(t)} \bigr)\\
&&\qquad\geq H(q_{i}-q_{i}',x) + \liminf_{t\to\infty} \bigl(v_{i}t - 2\eta
_{i}\sqrt{S_{i}t\log\log(2 S_{i}\eta_{i}^{2}t)} \bigr)\nonumber\\
&&\qquad= \infty\qquad\mbox{(almost surely)}.\nonumber
\end{eqnarray}
On the other hand, if $\lim_{t\to\infty}\rho_{i}(t)<\infty$, it is
trivial to obtain $G_{q_{i}-q_{i}'}(t)\to\infty$ by letting $t\to
\infty
$ in (\ref{eq:Gint}). Therefore, with $G_{q_{i}}(t)\geq
G_{q_{i}}(t)-G_{q'_{i}}(t)\to\infty$, we readily get $\lim_{t\to
\infty
}d_{\operatorname{KL}}(q_{i},X_{i}(t))=\infty$ (a.s.); and since
$G_{\alpha}(t) =
-\log
X_{i\alpha}(t)$ for all pure strategies $\alpha\in\mathcal{S}_{i}$, our
proof is complete.
\end{pf}

As in \cite{Im05}, we can now obtain the following estimate for the
lifespan of pure dominated strategies.
\begin{proposition}
\label{prop:timetolive}
Let $X(t)$ be a solution path of (\ref{eq:SRD}) with initial condition
$X(0) = x\in\operatorname{Int}(\Delta)$ and let $P_{x}$ denote its
law. Assume
further that the strategy $\alpha\in\mathcal{S}_{i}$ is dominated;
then, for
any $M>0$ and for $t$ large enough, we have
%
\begin{equation}
\label{eq:timetolive}
P_{x} \{X_{i\alpha}(t)<e^{-M} \} \geq\frac{1}{2} \operatorname
{erfc}\biggl(\frac{M -
h_{i}(x_{i}) - v_{i}t}{2\eta_{i}\sqrt{S_{i}t}} \biggr),
\end{equation}
where $S_{i}=|\mathcal{S}_{i}|$ is the number of strategies available to
player $i$, $\eta_{i} =\break \max\{|\eta_{i\beta}(y)|\dvtx y\in\Delta
, \beta\in
\mathcal{S}_{i}\}$ and the constants $v_{i}>0$ and $h_{i}(x_{i})$ do not
depend on $t$.
\end{proposition}
\begin{pf}
The proof is pretty straightforward and for the most part follows \cite
{Im05}. Surely enough, if $\alpha\prec p_{i}\in\Delta_{i}$ and we use
the same notation as in the proof of Proposition \ref{prop:dominated},
we have
%
\begin{eqnarray}
-\log X_{i\alpha}(t) &=& G_{\alpha}(t) \geq G_{\alpha}(t) -
G_{p_{i}}(t)\nonumber\\
&\geq& H(\alpha,x) - H(p_{i},x) + v_{i}t + \widetilde
{W}_{i}(\rho_{i}(t))\\
&=& h_{i}(x_{i}) + v_{i}t +\widetilde{W}_{i}(\rho_{i}(t)),\nonumber
\end{eqnarray}
where $v_{i} :=\min_{x_{-i}}\{u_{i}(x_{-i};p_{i}) -
u_{i}(x_{-i};\alpha
)\}>0$ and $h_{i}(x_{i}) := \log x_{i\alpha} -\break \sum_{\beta}p_{i\beta
}\log x_{i\beta}$. Then
%
\begin{eqnarray}
P_{x}\bigl(X_{i\alpha}(t) < e^{-M}\bigr) &\geq& P_{x} \{\widetilde{W}_{i}(\rho
_{i}(t))>M - h_{i}(x_{i}) - v_{i}t \}\nonumber\\[-8pt]\\[-8pt]
&=&\frac{1}{2}\operatorname{erfc}\biggl(\frac{M - h_{i}(x_{i}) -
v_{i}t}{\sqrt{2\rho
_{i}(t)}} \biggr)\nonumber
\end{eqnarray}
and, since the quadratic variation $\rho_{i}(t)$ is bounded above by $2
S_{i}\eta_{i}^{2}t$ (\ref{eq:quvar}), the estimate (\ref
{eq:timetolive}) holds for all sufficiently large $t$ [i.e., such that
$M<h_{i}(x_{i})+v_{i}t$].
\end{pf}

Some remarks are now in order: first and foremost, our results should
be contrasted to those of Cabrales \cite{Ca00} and Imhof \cite{Im05}
where dominated strategies die out only if the noise coefficients
(shocks) $\eta_{i\alpha}$ satisfy certain tameness conditions. The
origin of this notable difference is the form of the replicator
equation (\ref{eq:SRD}) and, in particular, the extra terms that are
propagated there by exponential learning and which are absent from the
aggregate shocks dynamics (\ref{eq:ASRD}). As can be seen from the
derivations in Proposition \ref{prop:dominated}, these terms are
precisely the ones that allow players to pick up on the true payoffs
$u_{i\alpha}$ instead of the modified ones $\widetilde{u}_{i\alpha} =
u_{i\alpha} + \frac{1}{2}\eta_{i\alpha}^{2}$ that come up in \cite
{Im05,Im09} (and, indirectly, in \cite{Ca00} as well).

Secondly, it turns out that the way that the noise coefficients $\eta
_{i\beta}$ depend on the profile $x\in\Delta$ is not really crucial:
as long as $\eta_{i\beta}(x)$ is continuous (or essentially bounded),
our arguments are not affected. The only way in which a specific
dependence influences the extinction of dominated strategies is seen in
Proposition \ref{prop:timetolive}: a sharper estimate of the quadratic
variation of $\int_{0}^{t} \eta_{i\beta}(X(s)) \,ds$ could conceivably
yield a more accurate estimate for the cumulative distribution function
of (\ref{eq:timetolive}).

Finally, it is only natural to ask if Proposition \ref{prop:dominated}
can be extended to strategies that are only \textit{iteratively}
dominated. As it turns out, this is indeed the case.
\begin{theorem}
\label{thm:rational}
Let $X(t)$ be a solution path of (\ref{eq:SRD}) starting at $X(0) =
x\in
\operatorname{Int}(\Delta)$. Then, if $q_{i}\in\Delta_{i}$ is
iteratively dominated,
%
\begin{equation}
\lim_{t\to\infty} d_{\operatorname{KL}}(q_{i},X_{i}(t)) = \infty
\qquad\mbox{almost surely,}
\end{equation}
that is, \textit{only rationally admissible strategies survive in the
long run.}
\end{theorem}
\begin{pf}
As in the deterministic case \cite{SZ92}, the main idea is that the
solution path $X(t)$ gets progressively closer to the faces of
$\Delta
$ that are spanned by the pure strategies which have not yet been
eliminated. Following \cite{Ca00}, we will prove this by induction on
the rounds of elimination of dominated strategies; Proposition \ref
{prop:dominated} is simply the case $n=1$.

To wit, let $A_{i}\subseteq\Delta_{i}$, $A_{-i}\subseteq\Delta
_{-i}$ and denote by $\operatorname{Adm}(A_{i},A_{-i})$ the set of strategies
$q_{i}\in A_{i}$ that are admissible (i.e., not dominated) with respect
to any strategy $q_{-i}\in A_{-i}$. So, if we start with $\mathcal
{A}_{i}^{0} =
\Delta_{i}$ and $\mathcal{A}_{-i}^{0} = \prod_{j\neq i}\mathcal
{A}_{j}^{0}$, we may
define inductively the set of strategies that remain admissible after
$n$ elimination rounds by $\mathcal{A}_{i}^{n} := \operatorname
{Adm}(\mathcal{A}_{i}^{n-1},\mathcal{A}
_{-i}^{n-1})$ where $\mathcal{A}_{i}^{n-1} := \prod_{j\neq i}
\mathcal{A}_{j}^{n-1}$;
similarly, the pure strategies that have survived after $n$ such rounds
will be denoted by $\mathcal{S}_{i}^{n}:= \mathcal{S}_{i}\cap
\mathcal{A}_{i}^{n}$.
Clearly, this sequence forms a descending chain $\mathcal
{A}_{i}^{0}\supseteq
\mathcal{A}
_{i}^{1}\supseteq\cdots$ and the set $\mathcal{A}_{i}^{\infty
}:=\bigcap
_{0}^{\infty}\mathcal{A}_{i}^{n}$ will consist precisely of the
strategies of
player $i$ that are rationally admissible.

Assume then that the cross entropy $G_{q_{i}}(t) = H(q_{i},X_{i}(t)) =
-\sum_{\alpha} q_{i\alpha}\times\break \log X_{i\alpha}(t)$ diverges as $t\to
\infty$ for all strategies $q_{i}\notin\mathcal{A}_{i}^{k}$ that die
out within
the first $k$ rounds; in particular, if $\alpha\notin\mathcal{S}_{i}^{k}$
this implies that $X_{i\alpha}(t)\to0$ as $t\to\infty$. We will show
that the same is true if $q_{i}$ survives for $k$ rounds but is
eliminated in the subsequent one.

Indeed, if $q_{i}\in\mathcal{A}_{i}^{k}$ but $q_{i}\notin\mathcal
{A}_{i}^{k+1}$, there
will exist some $q_{i}'\in\mathcal{A}_{i}^{k+1}$ such that
%
\begin{equation}
u_{i}(x_{-i};q_{i}') > u_{i}(x_{-i};q_{i})\qquad\mbox{for all }x_{-i}\in
\mathcal{A}
_{-i}^{k}.
\end{equation}
Now, note that any $x_{-i}\in\Delta_{-i}$ can be decomposed as
$x_{-i} = x_{-i}^{\mathrm{adm}} + x_{-i}^{\mathrm{dom}}$ where
$x_{-i}^{\mathrm{adm}}$ is the ``admissible'' part of $x_{-i}$, that is,
the projection of $x_{-i}$ on the subspace spanned by the surviving
vertices $\mathcal{S}^{k}_{-i} = \prod_{j\neq i}\mathcal
{S}_{i}^{k}$. Hence, if
$v_{i} = \min\{u_{i}(\alpha_{-i};q'_{i}) -
u_{i}(\alpha_{-i};q_{i})\dvtx
\alpha_{-i}\in\mathcal{S}_{-i}^{k}\}$, we will have $v_{i}>0$ and,
by linearity,
%
\begin{equation}
u_{i}(x_{-i}^{\mathrm{adm}};q_{i}') - u_{i}(x_{-i}^{\mathrm
{adm}};q_{i})\geq v_{i}>0\qquad\mbox{for all }x_{-i}\in\Delta_{-i}.
\end{equation}
Moreover, by the induction hypothesis, we also have $X_{-i}^{\mathrm
{dom}}(t)\to0$ as $t\to\infty$. Thus, there exists some $t_{0}$ such that
%
\begin{equation}
|u_{i}(X_{-i}^{\mathrm{dom}}(t), q_{i}') - u_{i}(X_{-i}^{\mathrm
{dom}}(t), q_{i}) |<v_{i}/2
\end{equation}
for all $t\geq t_{0}$ [recall that $X_{-i}^{\mathrm{dom}}(t)$ is
spanned by already eliminated strategies].

Therefore, as in the proof of Proposition \ref{prop:dominated}, we
obtain for $t\geq t_{0}$
%
\begin{equation}\qquad
G_{q_{i}}(t) - G_{q_{i}'}(t) \geq M + \frac{1}{2}v_{i} t + \sum
_{\beta
} (q_{i\beta}'-q_{i\beta})\int_{0}^{t} \eta_{i\beta}(X(s))
\,dW_{i\beta}(s),
\end{equation}
where $M$ is a constant depending only on $t_{0}$. In this way, the
same reasoning as before gives $\lim_{t\to\infty}G_{q_{i}}(t)=\infty$
and the theorem follows.
\end{pf}

As a result, if there exists only one rationally admissible strategy,
we get the following.
\begin{corollary}
\label{cor:dominance}
Let $X(t)$ be an interior solution path of the replicator equation
(\ref
{eq:SRD}) for some dominance-solvable game $\mathfrak{G}$ and let
$x_{0}\in
\mathcal{S}$ be the (unique) strict equilibrium of $\mathfrak{G}$. Then
%
\begin{equation}
\lim_{t\to\infty} X(t) = x_{0}\qquad \mbox{almost surely},
\end{equation}
that is, \textit{players converge to the game's strict equilibrium}
(\textit{a.s.}).
\end{corollary}

In concluding this section, it is important to note that all our
results on the extinction of dominated strategies remain true in the
adjusted dynamics (\ref{eq:SLRD}) as well: this is just a matter of
rescaling. The only difference in using different learning rates
$\lambda_{i}$ comes about in Proposition \ref{prop:timetolive} where
the estimate (\ref{eq:timetolive}) becomes
%
\begin{equation}
P_{x} \{X_{i\alpha}(t)<e^{-M} \} \geq\frac{1}{2} \operatorname
{erfc}\biggl(\frac{M -
h_{i}(x_{i}) - \lambda_{i} v_{i}t}{2\lambda_{i}\eta_{i}\sqrt
{S_{i}t}} \biggr).
\end{equation}

As it stands, this is not a significant difference in itself because
the two estimates are asymptotically equal for large times.
Nonetheless, it is this very lack of contrast that clashes with the
deterministic setting where faster learning rates accelerate the
emergence of rationality. The reason for this gap is that an increased
learning rate $\lambda_{i}$ also carries a commensurate increase in the
noise coefficients $\eta_{i}$, and thus deflates the benefits of
accentuating payoff differences. In fact, as we shall see in the next
sections, the learning rates do not really allow players to learn any
faster as much as they help diminish their shortsightedness: by
effectively being lazy, it turns out that players are better able to
average out the noise.

\section{Congestion games: A suggestive digression}
\label{sec:congestion}

Having established that irrational choices die out in the long run, we
turn now to the question of whether equilibrial play is stable in the
stochastic replicator dynamics of exponential learning. However, before
tackling this issue in complete generality, it will be quite
illustrative to pay a visit to the class of congestion games where the
presence of a potential simplifies things considerably. In this way,
the results we obtain here should be considered as a motivating
precursor to the general case analyzed in Section \ref{sec:equilibrium}.

\subsection{Congestion games}
To begin with, it is easy to see that the potential $V$ of Definition
\ref{def:potential} is a Lyapunov function for the deterministic
replicator dynamics. Indeed, assume that player $i\in\mathcal{N}$ is learning
at a rate $\lambda_{i}>0$ and let $x(t)$ be a solution path of the
rate-adjusted dynamics (\ref{eq:LRD}). Then, a simple differentiation
of $V(x(t))$ gives
%
\begin{eqnarray}
\frac{dV}{dt} &=& \sum_{i,\alpha}\frac{\partial V}{\partial
x_{i\alpha}}\,
\frac{d
x_{i\alpha}}{dt}
=-\sum_{i,\alpha} u_{i\alpha}(x) \lambda_{i} x_{i\alpha}
\bigl(u_{i\alpha
}(x) - u_{i}(x) \bigr)\nonumber\\[-8pt]\\[-8pt]
&=& -\sum_{i} \lambda_{i} \biggl(\sum_{\alpha} x_{i\alpha}
u^{2}_{i\alpha
}(x) - u_{i}^{2}(x) \biggr)\leq0,\nonumber
\end{eqnarray}
the last step following from Jensen's inequality---recall that $\frac
{\partial V}{\partial x_{i\alpha}} = - u_{i\alpha}(x)$ on account of
(\ref{eq:potential}) and also that $u_{i}(x) = \sum_{\alpha}
x_{i\alpha
} u_{i\alpha}(x)$. In particular, this implies that the trajectories
$x(t)$ are attracted to the local minima of $V$, and since these minima
coincide with the strict equilibria of the game, we painlessly infer
that strict equilibrial play is asymptotically stable in (\ref{eq:LRD})---as
mentioned before, we plead guilty to a slight abuse of
terminology in assuming that all equilibria in pure strategies are also strict.

It is therefore reasonable to ask whether similar conclusions can be
drawn in the noisy setting of (\ref{eq:SLRD}). Mirroring the
deterministic case, a promising way to go about this question is to
consider again the potential function $V$ of the game and try to show
that it is stochastically Lyapunov in the sense of Definition \ref
{def:Lyapunov}. Indeed, if $q_{0}=(e_{1,0},\ldots, e_{N,0})\in
\Delta$
is a local minimum of $V$ (and hence, a strict equilibrium of the
underlying game), we may assume without loss of generality that
$V(q_{0})=0$ so that $V(x)>0$ in a neighborhood of $q_{0}$. We are thus
left to examine the negativity condition of Definition \ref
{def:Lyapunov}, that is, whether there exists some $k>0$ such that $L
V(x)\leq-kV(x)$ for all $x$ sufficiently close to $q_{0}$.

To that end, recall that $\frac{\partial V}{\partial x_{i\alpha}} =
-u_{i\alpha}$
and that $\frac{\partial^{2}V}{\partial x_{i\alpha}^{2}} =0$. Then,
the generator
$L$ of the rate-adjusted dynamics (\ref{eq:SLRD}) applied to $V$ produces
%
\begin{eqnarray}
\label{eq:genV}\qquad
LV(x) &=&-\sum_{i,\alpha} \lambda_{i} x_{i\alpha} u_{i\alpha}(x)
\bigl(u_{i\alpha}(x) - u_{i}(x) \bigr)\nonumber\\[-8pt]\\[-8pt]
&&{}-\sum_{i,\alpha}
\frac{\lambda_{i}^{2}}{2} x_{i\alpha} u_{i\alpha}(x)
\biggl(\eta_{i\alpha}^{2}(1-2x_{i\alpha}) - \sum_{\beta}\eta_{i\beta}^{2}
x_{i\beta}(1-2x_{i\beta}) \biggr),\nonumber
\end{eqnarray}
where, for simplicity, we have assumed that the noise coefficients
$\eta
_{i\alpha}$ are constant.

We will study (\ref{eq:genV}) term by term by considering the perturbed
strategies $x_{i} = (1-\varepsilon_{i}) e_{i,0} + \varepsilon_{i}
y_{i}$ where
$y_{i}$ belongs to the face of $\Delta_{i}$ that lies opposite to
$e_{i,0}$ (i.e., $y_{i\mu}\geq0$, $\mu=1,2,\ldots$ and $\sum_{\mu}
y_{i\mu}=1$) and $\varepsilon_{i}>0$ measures the distance of player
$i$ from
$e_{i,0}$. In this way, we get
%
\begin{eqnarray}
u_{i}(x) &=& \sum_{\alpha} x_{i\alpha} u_{i\alpha} (x)
= (1-\varepsilon_{i}) u_{i,0}(x) + \varepsilon_{i} \sum_{\mu}
y_{i\mu} u_{i\mu
}(x)\nonumber\\
&=& u_{i,0}(x)
+\varepsilon_{i} \sum_{\mu} y_{i\mu} [u_{i\mu}(x) - u_{i,0}(x)]\\
&=& u_{i,0}(x)
-\varepsilon_{i} \sum_{\mu} y_{i\mu} \Delta u_{i\mu}
+ \mathcal O(\varepsilon_{i}^{2} ),\nonumber
\end{eqnarray}
where $\Delta u_{i\mu} = u_{i,0}(q_{0}) - u_{i\mu}(q_{0})>0$. Then, by
going back to (\ref{eq:genV}), we obtain
%
\begin{eqnarray}
\label{eq:paystep}
&&\sum_{\alpha} x_{i\alpha} u_{i\alpha}(x) [u_{i\alpha}(x) - u_{i}(x)
]\nonumber\\
&&\qquad= (1-\varepsilon_{i}) u_{i,0}(x) [u_{i,0}(x) - u_{i}(x) ]\nonumber\\
&&\qquad\quad{}
+ \varepsilon_{i} \sum_{\mu} y_{i\mu} u_{i\mu} (x) [u_{i\mu}(x) - u_{i}(x)
]\nonumber\\
&&\qquad= (1-\varepsilon_{i}) u_{i,0}(x) \cdot\varepsilon_{i} \sum_{\mu}
y_{i\mu
}\Delta
u_{i\mu}\\
&&\qquad\quad{}- \varepsilon_{i} \sum_{\mu} y_{i\mu} u_{i\mu}(q_{0}) \Delta
u_{i\mu}
+ \mathcal O(\varepsilon_{i}^{2} )\nonumber\\
&&\qquad= \varepsilon_{i} \sum_{\mu} y_{i\mu} u_{i,0}(q_{0}) \Delta
u_{i\mu}
- \varepsilon_{i} \sum_{\mu} y_{i\mu} u_{i\mu}(q_{0}) \Delta
u_{i\mu}
+ \mathcal O(\varepsilon_{i}^{2} )\nonumber\\
&&\qquad= \varepsilon_{i} \sum_{\mu} y_{i\mu} (\Delta u_{i\mu} )^{2}
+ \mathcal O(\varepsilon_{i}^{2} ).\nonumber
\end{eqnarray}
As for the second term of (\ref{eq:genV}), some easy algebra reveals that
%
\begin{eqnarray}
&&\eta_{i,0}^{2}(1-2x_{i,0})
- \sum_{\beta}\eta_{i\beta}^{2} x_{i\beta}(1-2x_{i\beta
})\nonumber\\
&&\qquad= -\eta_{i,0}^{2} (1-2\varepsilon_{i})
- \eta_{i,0}^{2}(1-\varepsilon_{i}) - \varepsilon_{i}\sum_{\mu}
\eta_{i\mu}^{2}
y_{i\mu}\nonumber\\[-8pt]\\[-8pt]
&&\qquad\quad{} + 2(1-\varepsilon_{i})^{2} \eta_{i,0}^{2}
+ 2\varepsilon_{i}^{2}\sum_{\mu} \eta_{i\mu}^{2} y_{i\mu
}^{2}\nonumber\\
&&\qquad= -\varepsilon_{i} \biggl(\eta_{i,0}^{2} + \sum_{\mu} y_{i\mu} \eta
_{i\mu
}^{2} \biggr)
+ \mathcal O(\varepsilon_{i}^{2} )\nonumber
\end{eqnarray}
and, after a (somewhat painful) series of calculations, we get
%
\begin{eqnarray}
\label{eq:noisestep}
&&\sum_{\alpha}
x_{i\alpha} u_{i\alpha}(x) \biggl(
\eta_{i\alpha}^{2}(1-2x_{i\alpha}) - \sum_{\beta}\eta_{i\beta}^{2}
x_{i\beta}(1-2x_{i\beta})
\biggr)\nonumber\\
&&\qquad= (1-\varepsilon_{i})
u_{i,0}(x) \biggl(
\eta_{i,0}^{2}(1-2x_{i,0}) - \sum_{\beta}\eta_{i\beta}^{2}
x_{i\beta
}(1-2x_{i\beta})
\biggr)\nonumber\\
&&\qquad\quad{}+ \varepsilon_{i}\sum_{\mu}y_{i\mu} \biggl(
\eta_{i\mu}^{2}(1-2x_{i\mu}) - \sum_{\beta}\eta_{i\beta}^{2}
x_{i\beta}(1-2x_{i\beta})
\biggr)\nonumber\\[-8pt]\\[-8pt]
&&\qquad= -\varepsilon_{i} u_{i,0}(q_{0}) \biggl(
\eta_{i,0}^{2} + \sum_{\mu}y_{i\mu}\eta_{i\mu}^{2}
\biggr)\nonumber\\
&&\qquad\quad{}+ \varepsilon_{i}\sum_{\mu} y_{i\mu} u_{i\mu}(q_{0}) (\eta_{i\mu
}^{2} +
\eta_{i,0}^{2})
+ \mathcal O(\varepsilon_{i}^{2} )\nonumber\\
&&\qquad= -\varepsilon_{i}\sum_{\mu} y_{i\mu} \Delta u_{i\mu} (\eta
_{i\mu
}^{2} +
\eta_{i,0}^{2} )
+ \mathcal O(\varepsilon_{i}^{2} ).\nonumber
\end{eqnarray}
Finally, if we assume without loss of generality that $V(q_{0})=0$ and
set $\xi= x - q_{0}$ (i.e., $\xi_{i,0} = -\varepsilon_{i}$ and $\xi
_{i\mu
} =
\varepsilon_{i} y_{i\mu}$ for all $i\in\mathcal{N}$, $\mu\in
\mathcal{S}_{i}
\setminus\{
0\}$), we readily get
%
\begin{eqnarray}
\label{eq:potstep}
V(x) &=& \sum_{i,\alpha} \frac{\partial V}{\partial x_{i\alpha}} \xi
_{i\alpha}
+ \mathcal O(\xi^{2} )\nonumber\\
&=& -\sum_{i,\alpha}\frac{\partial u_{i}}{\partial x_{i\alpha}}
\bigg|_{q_{0}} \xi
_{i\alpha}
+ \mathcal O(\varepsilon^{2} )\nonumber\\[-8pt]\\[-8pt]
&=& -\sum_{i,\alpha} u_{i\alpha}(q_{0}) \xi_{i\alpha}
+ \mathcal O(\varepsilon^{2} )\nonumber\\
&=& \sum_{i}\varepsilon_{i}\sum_{\mu} y_{i\mu}\Delta u_{i\mu}
+ \mathcal O(\varepsilon^{2} ),\nonumber
\end{eqnarray}
where $\varepsilon^{2} = \sum_{i} \varepsilon_{i}^{2}$. Therefore,
by combining
(\ref
{eq:paystep}), (\ref{eq:noisestep}) and (\ref{eq:potstep}), the
negativity condition $LV(x) \leq-k V(x)$ becomes
%
\begin{eqnarray}
&&\sum_{i}\lambda_{i}\varepsilon_{i}
\sum_{\mu} y_{i\mu}\Delta u_{i\mu} \biggl[
\Delta u_{i\mu}
- \frac{\lambda_{i}}{2} (\eta_{i\mu}^{2} + \eta_{i,0}^{2} )
\biggr]\nonumber\\[-8pt]\\[-8pt]
&&\qquad\geq k\sum_{i}\varepsilon_{i} \sum_{\mu} y_{i\mu}\Delta u_{i\mu} +
\mathcal O
(\varepsilon^{2}).\nonumber
\end{eqnarray}
Hence, if $\Delta u_{i\mu} > \frac{\lambda_{i}}{2} (\eta_{i\mu
}^{2} +
\eta_{i,0}^{2})$ for all $\mu\in\mathcal{S}_{i} \setminus\{0\}$,
this last
inequality will be satisfied for some $k>0$ whenever $\varepsilon$ is small
enough. Essentially, this proves the following.
\begin{proposition}
\label{prop:potential}
Let $q=(\alpha_{1},\ldots,\alpha_{N})$ be a strict equilibrium of a
congestion game $\mathfrak{G}$ with potential function $V$ and assume that
$V(q)=0$. Assume further that the learning rates $\lambda_{i}$ are
sufficiently small so that, for all $\mu\in\mathcal{S}_{i} \setminus
\{
\alpha
_{i}\}$ and all $i\in\mathcal{N}$,
%
\begin{equation}
\label{eq:payvsnoise}
V(q_{-i},\mu) > \frac{\lambda_{i}}{2}(\eta_{i\mu}^{2} +\eta_{i,0}^{2}).
\end{equation}
Then $q$ is stochastically asymptotically stable in the rate-adjusted
dynamics (\ref{eq:SLRD}).
\end{proposition}

We thus see that no matter how loud the noise $\eta_{i}$ might be,
stochastic stability is always guaranteed if the players choose a
learning rate that is slow enough as to allow them to average out the
noise (i.e., $\lambda_{i}<\Delta V_{i}/\eta_{i}^{2}$). Of course, it
can be argued here that it is highly unrealistic to expect players to
be able to estimate the amount of Nature's interference and choose a
suitably small rate $\lambda_{i}$. On top of that, the very form of the
condition (\ref{eq:payvsnoise}) is strongly reminiscent of the
``modified'' game of \mbox{\cite{Im05,Im09}}, a similarity which seems to
contradict our statement that exponential learning favors rational
reactions in the \textit{original} game. The catch here is that condition
(\ref{eq:payvsnoise}) is only \textit{sufficient} and Proposition
\ref{prop:potential} merely highlights the role of a potential function in
a stochastic environment. As we shall see in Section \ref
{sec:equilibrium}, nothing stands in the way of choosing a different
Lyapunov candidate and dropping requirement (\ref{eq:payvsnoise}) altogether.

\subsection{The dyadic case}

To gain some further intuition into why the condition (\ref
{eq:payvsnoise}) is redundant, it will be particularly helpful to
examine the case where players compete for the resources of only two
facilities (i.e., $\mathcal{S}_{i} = \{0,1\}$ for all $i\in\mathcal
{N}$) and try
to learn the game with the help of the uniform replicator equation
(\ref
{eq:SRD}). This is the natural setting for the El Farol bar problem
\cite{Ar94} and the ensuing minority game \cite{MC00} where players
choose to ``buy'' or ``sell'' and are rewarded when they are in the
minority--buyers in a sellers' market or sellers in an abundance of buyers.

As has been shown in \cite{Mi96}, such games always possess strict
equilibria, even when players have distinct payoff functions. So,
by
relabeling indices if necessary, let us assume that $q_{0} =
(e_{1,0},\ldots, e_{N,0})$ is such a strict equilibrium and set \mbox{$x_{i}
\equiv x_{i,0}$}. Then, the generator of the replicator equation (\ref
{eq:SRD}) takes the form
%
\begin{eqnarray}
\label{eq:2Dgen}
L&=&\sum_{i} x_{i}(1-x_{i}) \biggl[\Delta u_{i}(x)+ \frac
{1}{2}(1-2x_{i})\eta_{i}^{2}(x) \biggr]\, \frac{\partial}{\partial
x_{i}}\nonumber\\[-8pt]\\[-8pt]
&&{}+ \frac{1}{2}\sum_{i}x_{i}^{2}(1-x_{i})^{2} \eta_{i}^{2}(x)\, \frac
{\partial^{2}}{\partial x_{i}^{2}},\nonumber
\end{eqnarray}
where now $\Delta u_{i} \equiv u_{i,0} - u_{i,1}$ and $\eta_{i}^{2} =
\eta_{i,0}^{2} + \eta_{i,1}^{2}$.

It thus appears particularly appealing to introduce a new set of
variables $y_{i}$ such that $\frac{\partial}{\partial y_{i}} =
x_{i}(1-x_{i})\,\frac
{\partial}{\partial x_{i}}$; this is just the ``logit''
transformation: $y_{i} =
\operatorname{logit}x_{i} \equiv\log\frac{x_{i}}{1-x_{i}}$. In
these new
variables, (\ref{eq:2Dgen}) assumes the astoundingly suggestive guise
%
\begin{equation}
\label{eq:2Dygen}
L= \sum_{i} \biggl(\Delta u_{i}\, \frac{\partial}{\partial y_{i}} + \frac{1}{2}
\eta_{i}^{2}\, \frac{\partial^{2}}{\partial y_{i}^{2}} \biggr),
\end{equation}
which reveals that the noise coefficients can be effectively decoupled
from the payoffs. We can then take advantage of this by letting $L$
act on the function $f(y) = \sum_{i} e^{-a_{i}y_{i}}$ ($a_{i}>0$):
%
\begin{equation}
Lf(y) = -\sum_{i} a_{i} \biggl(\Delta u_{i} - \frac{1}{2}a_{i}\eta
_{i}^{2} \biggr) e^{-a_{i}y_{i}}.
\end{equation}
Hence, if $a_{i}$ is chosen small enough so that $\Delta u_{i} - \frac
{1}{2}a_{i}\eta_{i}^{2}\geq m_{i}>0$ for all sufficiently large $y_{i}$
[recall that $\Delta u_{i}(q_{0})>0$ since $q_{0}$ is a strict
equilibrium], we get
%
\begin{equation}
Lf(y) \leq-\sum_{i} a_{i}m_{i} e^{-a_{i}y_{i}} \leq- k f(y),
\end{equation}
where $k = \min_{i}\{a_{i} m_{i}\}>0$. And since $f$ is strictly
positive for $y_{i,0}>0$ and only vanishes as $y\to\infty$ (i.e., at
the equilbrium $q_{0}$), a trivial modification of the stochastic
Lyapunov method (see, e.g., pages 314 and 315 of \cite{GS71}) yields
the following.
\begin{proposition}
\label{prop:minority}
The strict equilibria of minority games are stochastically
asymptotically stable in the uniform replicator equation (\ref{eq:SRD}).
\end{proposition}
\begin{remark}
It is trivial to see that strict equilibria of minority games will also
be stable in the rate-adjusted\vspace*{1pt} dynamics (\ref{eq:SLRD}): in that case
we simply need to choose $a_{i}$ such that $\Delta u_{i} -\frac
{1}{2}a_{i}\lambda_{i} \eta_{i}^{2}\geq m_{i} >0$.
\end{remark}
\begin{remark}
A closer inspection of the calculations leading to Proposition~\ref
{prop:minority} reveals that nothing hinges on the minority mechanism
per se: it is (\ref{eq:2Dygen}) that is crucial to our analysis and
$L$ takes this form whenever the underlying game is a \textit{dyadic}
one (i.e., $|\mathcal{S}_{i}| = 2$ for all $i\in\mathcal{N}$). In
other words,
Proposition \ref{prop:minority} also holds for all games with 2
strategies and should thus be seen as a significant extension of
Proposition \ref{prop:potential}.
\end{remark}
\begin{proposition}
\label{prop:dyadic}
The strict equilibria of dyadic games are stochastically a\-sym\-
ptotically stable in the replicator dynamics (\ref{eq:SRD}), (\ref
{eq:SLRD}) of exponential learning.
\end{proposition}

\section{Stability of equilibrial play}
\label{sec:equilibrium}

In deterministic environments, the ``folk theorem'' of evolutionary
game theory provides some pretty strong ties between equilibrial play
and stability: strict equilibria are asymptotically stable in the
multi-population replicator dynamics (\ref{eq:RD}) \cite{We95}. In our
stochastic setting, we have already seen that this is always true in
two important classes of games: those that can be solved by iterated
elimination of dominated strategies (Corollary \ref{cor:dominance}) and
dyadic ones (Proposition \ref{prop:dyadic}).

Although interesting in themselves, these results clearly fall short of
adding up to a decent analogue of the folk theorem for stochastically
perturbed games. Nevertheless, they are quite strong omens in that
direction and such expectations are vindicated in the following.
\begin{theorem}
\label{thm:stability}
The strict equilibria of a game $\mathfrak{G}$ are stochastically
asymptotically stable in the replicator dynamics (\ref{eq:SRD}), (\ref
{eq:SLRD}) of exponential learning.
\end{theorem}

Before proving Theorem \ref{thm:stability}, we should first take a
slight detour in order to properly highlight some of the issues at
hand. On that account, assume again that the profile
$q_{0}=(e_{1,0},\ldots,e_{N,0})$ is a strict equilibrium of $\mathfrak{G}$.
Then, if $q_{0}$ is to be stochastically stable, say in the uniform
dynamics (\ref{eq:SRD}), one would expect the strategy scores $U_{i,0}$
of player $i$ to grow much faster than the scores $U_{i\mu},\mu\in
\mathcal{S}
_{i} \setminus\{0\}$ of his other strategies. This is captured
remarkably well by the ``adjusted'' scores
%
\begin{subequation}
\label{eq:ratescore}
\begin{eqnarray}
Z_{i,0} &=& \lambda_{i} U_{i,0} - \log\biggl(\sum_{\mu} e^{\lambda
_{i}U_{i\mu}} \biggr),\\
Z_{i\mu} &=& \lambda_{i}(U_{i\mu} - U_{i,0}),
\end{eqnarray}
\end{subequation}
where $\lambda_{i}>0$ is a sensitivity parameter akin (but not
identical) to the learning rates of (\ref{eq:SLRD}) (the
choice of common notation is fairly premeditated though).

Clearly, whenever $Z_{i,0}$ is large, $U_{i,0}$ will be much greater
than any other score $U_{i\mu}$ and hence, the strategy $0\in\mathcal{S}
_{i}$ will be employed by player $i$ far more often. To see this in
more detail, it is convenient to introduce the variables
%
\begin{subequation}
\label{eq:Ydef}
\begin{eqnarray}
\label{eq:Y0}
Y_{i,0} &:=& e^{Z_{i,0}}
= \frac{e^{\lambda_{i}U_{i,0}}}{\sum_{\nu} e^{\lambda_{i} U_{i\nu
}}},\\
\label{eq:Ym}
Y_{i\mu} &:=& \frac{e^{Z_{i\mu}}}{\sum_{\nu} e^{Z_{i\nu}}}
= \frac{e^{\lambda_{i} U_{i\mu}}}{\sum_{\nu} e^{\lambda_{i}
U_{i\nu}}},
\end{eqnarray}
\end{subequation}
where $Y_{i,0}$ is a measure of how close $X_{i}$ is to $e_{i,0}\in
\Delta_{i}$ and $(Y_{i,1}, Y_{i,2},\ldots)\in\Delta^{S_{i}-1}$
is a
direction indicator; the two sets of coordinates are then related by
the transformation $Y_{i\alpha} = X_{i\alpha}^{\lambda_{i}}/\sum
_{\mu}
X_{i\mu}^{\lambda_{i}}$, $\alpha\in\mathcal{S}_{i}$, $\mu\in
\mathcal{S}_{i}
\setminus\{0\}$. Consequently, to show that the strict equilibrium
$q_{0}=(e_{1,0},\ldots, e_{N,0})$ is stochastically asymptotically
stable in the replicator equation (\ref{eq:SRD}), it will suffice to
show that $Y_{i,0}$ diverges to infinity as $t\to\infty$ with
arbitrarily high probability.

Our first step in this direction will be to derive an SDE for the
evolution of the $Y_{i\alpha}$ processes. To that end, It\^{o}'s lemma gives
%
\begin{eqnarray}
\label{eq:ItoY}
d Y_{i\alpha}
&=& \sum_{j,\beta} \frac{\partial Y_{i\alpha}}{\partial U_{j\beta}}
\,dU_{j\beta}
+ \frac{1}{2} \sum_{j,k}\sum_{\beta,\gamma}
\frac{\partial^{2} Y_{i\alpha}}{\partial U_{j\beta}\,\partial
U_{k\gamma}}
\,dU_{j\beta} \cdot dU_{k\gamma}\nonumber\\[-8pt]\\[-8pt]
&=& \sum_{\beta} \biggl(
u_{i\beta}\frac{\partial Y_{i\alpha}}{\partial U_{i\beta}}
+ \frac{1}{2} \eta_{i\beta}^{2} \,\frac{\partial^{2} Y_{i\alpha
}}{\partial
U_{i\beta}^{2}}
\biggr) \,dt
+ \sum_{\beta} \eta_{i\beta} \,\frac{\partial Y_{i\alpha}}{\partial
U_{i\beta}}
\,dW_{i\beta},\nonumber
\end{eqnarray}
where, after a simple differentiation of (\ref{eq:Y0}), we have
\begin{subequation}
\begin{eqnarray}
\label{eq:Y0U0}
\hspace*{-23.5pt}\frac{\partial Y_{i,0}}{\partial U_{i,0}}
&=& \lambda_{i} Y_{i,0},\qquad
\frac{\partial^{2} Y_{i,0}}{\partial U_{i,0}^{2}}
= \lambda_{i}^{2} Y_{i,0},\hspace*{18pt}\nonumber\\[-34pt]
\end{eqnarray}
\renewcommand{\theequation}{6.4a$'$}
\begin{eqnarray}
\hspace*{23.5pt}\frac{\partial Y_{i,0}}{\partial U_{i\nu}}
&=& -\lambda_{i} Y_{i,0} Y_{i\nu},\qquad
\frac{\partial^{2} Y_{i,0}}{\partial U_{i\nu}^{2}}
= -\lambda_{i}^{2} Y_{i,0} Y_{i\nu}(1-2 Y_{i\nu})\hspace*{-23.5pt}
\end{eqnarray}
\end{subequation}
and, similarly, from (\ref{eq:Ym})
{\renewcommand{\theequation}{6.4b}
\begin{eqnarray}
\label{eq:YmU0}
\hspace*{8.5pt}\frac{\partial Y_{i\mu}}{\partial U_{i,0}}
&=& 0,\qquad
\frac{\partial^{2} Y_{i\mu}}{\partial U_{i,0}^{2}}
= 0,\hspace*{47.5pt}\nonumber\\[-34pt]
\end{eqnarray}
\renewcommand{\theequation}{6.4b$'$}
\begin{eqnarray}
\frac{\partial Y_{i\mu}}{\partial U_{i\nu}}
&=& \lambda_{i} Y_{i\mu}(\delta_{\mu\nu} -
Y_{i\nu}),\nonumber\\[-8pt]\\[-8pt]
\frac{\partial^{2} Y_{i\mu}}{\partial U_{i\nu}^{2}}
&=& \lambda_{i}^{2} Y_{i\mu}(\delta_{\mu\nu} -
Y_{i\nu})(1-2Y_{i\nu}).\nonumber
\end{eqnarray}}

\noindent In this way, by plugging everything back into (\ref{eq:ItoY}) we
finally obtain
\setcounter{equation}{4}
\begin{subequation}
\label{eq:dY}
\begin{eqnarray}
\label{eq:dY0}\hspace*{37pt}
dY_{i,0} &=& \lambda_{i}Y_{i,0} \biggl[
u_{i,0} - \sum_{\mu} Y_{i\mu} u_{i\mu}
+ \frac{\lambda_{i}}{2}\eta_{i,0}^{2}
- \frac{\lambda_{i}}{2} \sum_{\mu}Y_{i\mu}(1-2Y_{i\mu})\eta
_{i\mu}^{2}
\biggr] \,dt\nonumber\\[-8pt]\\[-8pt]
&&{}+\lambda_{i} Y_{i,0} \biggl[
\eta_{i,0} \,dW_{i,0} -\sum_{\mu} \eta_{i\mu} Y_{i\mu} \,dW_{i\mu}
\biggr],\nonumber\\
\label{eq:dYm}
dY_{i\mu} &=& \lambda_{i}Y_{i\mu}
[u_{i\mu} - \sum_{\nu} u_{i\nu} Y_{i\nu} ] \,dt\nonumber\\
&&{}+\frac{\lambda_{i}^{2}}{2}Y_{i\mu}
\biggl[\eta_{i\mu}^{2}(1-2 Y_{i\mu}) -\sum_{\nu}\eta_{i\nu}^{2}
Y_{i\nu
}(1-2Y_{i\nu}) \biggr] \,dt\\
&&{}+\lambda_{i} Y_{i\mu} \biggl[\eta_{i\mu} \,dW_{i\mu} - \sum_{\nu}
\eta
_{i\nu} Y_{i\nu} \,dW_{i\nu} \biggr],\nonumber
\end{eqnarray}
\end{subequation}
where we have suppressed the arguments of $u_{i}$ and $\eta_{i}$ in
order to reduce notational clutter.

This last SDE is particularly revealing: roughly speaking, we see that
if $\lambda_{i}$ is chosen small enough, the deterministic term
$u_{i,0} - \sum_{\mu} Y_{i\mu} u_{i\mu}$ will dominate the rest (cf.
with the ``soft'' learning rates of Proposition \ref{prop:potential}).
And, since we know that strict equilibria are asymptotically stable in
the deterministic case, it is plausible to expect the SDE (\ref{eq:dY})
to behave in a similar fashion.
\begin{pf*}{Proof of Theorem \protect\ref{thm:stability}}
Tying in with our previous discussion, we will establish stochastic
asymptotic stability of strict equilibria in the dynamics (\ref
{eq:SRD}) by looking at the processes $Y_{i} = (Y_{i,0}, Y_{i,1},\ldots
) \in{\mathbb R}\times\Delta^{S_{i}-1}$ of (\ref{eq:Ydef}). In these
coordinates, we just need to show that for every $M_{i}>0, i\in
\mathcal{N}$
and any $\varepsilon>0$, there exist $Q_{i} >M_{i}$ such that if
$Y_{i,0}(0) >
Q_{i}$, then, with probability greater than $1-\varepsilon$, $\lim
_{t\to
\infty
}Y_{i,0}(t) = \infty$ and $Y_{i,0}(t) > M_{i}$ for all $t\geq0$. In the
spirit of the previous section, we will accomplish this with the help
of the stochastic Lyapunov method.

Our first task will be to calculate the generator of the diffusion
$Y=(Y_{1},\ldots, Y_{N})$, that is, the second order differential operator
%
\begin{equation}
L= \mathop{\sum_{i\in\mathcal{N}}}_{\alpha\in\mathcal{S}_{i}}
b_{i\alpha} (y) \,\frac{\partial}{\partial y_{i\alpha}}
+\frac{1}{2} \mathop{\sum_{i\in\mathcal{N}}}_{\alpha,\beta\in
\mathcal{S}_{i}}
(\sigma_{i}(y)\sigma_{i}^{T}(y) )_{\alpha\beta}\,
\frac{\partial^{2}}{\partial y_{i\alpha}\, \partial y_{i\beta}},
\end{equation}
where $b_{i}$ and $\sigma_{i}$ are the drift and diffusion coefficients
of the SDE (\ref{eq:dY}), respectively.
In particular, if we restrict our attention to sufficiently smooth
functions of the form $f(y) = \sum_{i\in\mathcal{N}}
f_{i}(y_{i,0})$, the
application of $L$ yields
%
\begin{eqnarray}
\label{eq:generator}
Lf(y) &=& \sum_{i\in\mathcal{N}} \lambda_{i} y_{i,0} \biggl[
u_{i,0} +\frac{\lambda_{i}}{2}\eta_{i,0}^{2}
\nonumber\\
&&\hspace*{46.2pt}{}
- \sum_{\mu} y_{i\mu} \biggl(
u_{i\mu} -\frac{\lambda_{i}}{2} (1-2y_{i\mu})\eta_{i\mu}^{2} \biggr)
\biggr]\,
\frac{\partial f_{i}}{\partial y_{i,0}}\\
&&{} +\frac{1}{2} \sum_{i\in\mathcal{N}} \lambda_{i}^{2} y_{i,0}^{2} \biggl[
\eta_{i,0}^{2} + \sum_{\mu} \eta_{i\mu}^{2} y_{i\mu}^{2} \biggr]\,
\frac{\partial^{2} f_{i}}{\partial^{2}y_{i,0}}.\nonumber
\end{eqnarray}

Therefore, let us consider the function $f(y) = \sum_{i}1/y_{i,0}$ for
$ y_{i,0}>0$. With $\frac{\partial f}{\partial y_{i,0}} =
-1/y_{i,0}^{2}$ and
$\frac{\partial^{2} f}{\partial y_{i,0}^{2}} = 2/y_{i,0}^{3}$, (\ref
{eq:generator}) becomes
%
\begin{eqnarray}
Lf(y) &=& -\sum_{i\in\mathcal{N}} \frac{\lambda_{i}}{y_{i,0}} \biggl[
u_{i,0} - \sum_{\mu} u_{i\mu} y_{i\mu}
-\frac{\lambda_{i}}{2} \eta_{i,0}^{2}\nonumber\\[-8pt]\\[-8pt]
&&\hspace*{54.1pt}{} - \frac{\lambda_{i}}{2}\sum_{\mu} y_{i\mu}(1-y_{i\mu})\eta
_{i\mu
}^{2} \biggr].\nonumber
\end{eqnarray}
However, since $q_{0}=(e_{1,0},\ldots, e_{N,0})$ has been assumed to be
a strict Nash equilibrium of $\mathfrak{G}$, we will have
$u_{i,0}(q_{0})>u_{i\mu}(q_{0})$ for all $\mu\in\mathcal{S}_{i}
\setminus\{
0\}$. Then, by continuity, there exists some positive constant
$v_{i}>0$ with $u_{i,0} -\sum_{\mu}u_{i\mu}y_{i\mu} \geq v_{i}>0$
whenever $y_{i,0}$ is large enough (recall that $\sum_{\mu}y_{i\mu
}=1$). So, if we set $\eta_{i} = \max\{|\eta_{i\beta}(x)|\dvtx x\in
\Delta,
\beta\in\mathcal{S}_{i}\}$ and pick positive $\lambda_{i}$ with
$\lambda_{i}
<v_{i}/\eta_{i}^{2}$, we get
%
\begin{equation}
Lf(y) \leq-\sum_{i\in\mathcal{N}} \frac{\lambda_{i} v_{i}}{2}
\frac
{1}{y_{i,0}} \leq-\frac{1}{2}\min_{i}\{\lambda_{i} v_{i}\} f(y)
\end{equation}
for all sufficiently large $y_{i,0}$. Moreover, $f$ is strictly
positive for $y_{i,0}>0$ and vanishes only as $y_{i,0}\to\infty$.
Hence, as in the proof of Proposition \ref{prop:minority}, our claim
follows on account of $f$ being a (local) stochastic Lyapunov function.

Finally, in the case of the rate-adjusted replicator dynamics (\ref
{eq:SLRD}), the proof is similar and only entails a rescaling of the
parameters $\lambda_{i}$.
\end{pf*}
\begin{remark}
If we trace our steps back to the coordinates $X_{i\alpha}$, our
Lyapunov candidate takes the form $f(x) = \sum_{i} (x_{i,0}^{-\lambda
_{i}}\sum_{\mu} x_{i\mu}^{\lambda_{i}} )$. It thus begs to be compared
to the Lyapunov function $\sum_{\mu} x_{\mu}^{\lambda}$ employed by
Imhof and Hofbauer in \cite{Im09} to derive a conditional version of
Theorem \ref{thm:stability} in the evolutionary setting. As it turns
out, the obvious extension $f(x) = \sum_{i}\sum_{\mu} x_{i\mu
}^{\lambda
_{i}}$ works in our case as well, but the calculations are much more
cumbersome and they are also shorn of their ties to the adjusted scores
(\ref{eq:ratescore}).
\end{remark}
\begin{remark}
We should not neglect to highlight the dual role that the learning
rates $\lambda_{i}$ play in our analysis. In the logistic learning
model (\ref{eq:ratelogit}), they measure the players' convictions and
how strongly they react to a given stimulus (the scores $U_{i\alpha}$);
in this role, they are fixed at the outset of the game and form an
intrinsic part of the replicator dynamics (\ref{eq:SLRD}). On the other
hand, they also make a virtual appearance as free temperature
parameters in the adjusted scores (\ref{eq:ratescore}), to be softened
until we get the desired result. For this reason, even though Theorem
\ref{thm:stability} remains true for any\vspace*{1pt} choice of learning rates, the
function $f(x) = \sum_{i} x_{i,0}^{-\lambda_{i}}\sum_{\mu} x_{i\mu
}^{\lambda_{i}}$ is Lyapunov only if the sensitivity parameters
$\lambda
_{i}$ are small enough. It might thus seem unfortunate that we chose
the same notation in both cases, but we feel that our decision is
justified by the intimate relation of the two parameters.
\end{remark}

\section{Discussion}
\label{sec:discussion}

Our aim in this last section will be to discuss a number of important
issues that we have not been able to address thoroughly in the rest of
the paper; truth be told, a good part of this discussion can be seen as
a roadmap for future research.

\subsection*{Ties with evolutionary game theory}

In single-population evolutionary models, an evolutionarily stable
strategy (ESS) is a strategy which is robust against invasion by mutant
phenotypes \cite{MS74}. Strategies of this kind can be considered as a
stepping stone between mixed and strict equilibria and they are of such
significance that it makes one wonder why they have not been included
in our analysis.

The reason for this omission is pretty simple: even the weakest
evolutionary criteria in multi-population models tend to reject all
strategies which are not strict Nash equilibria \cite{We95}. Therefore,
since our learning model (\ref{eq:RD}) corresponds exactly to the
multi-population environment (\ref{eq:ERD}), we lose nothing by
concentrating our analysis only on the strict equilibria of the game.
If anything, this equivalence between ESS and strict equilibria in
multi-population settings further highlights the importance of the latter.

However, this also brings out the gulf between the single-population
setting and our own, even when we restrict ourselves to 2-player games
(which are the norm in single-population models). Indeed, the
single-population version of the dynamics (\ref{eq:ASRD}) is:
%
\begin{eqnarray}
\label{eq:1ASRD}
dX_{\alpha} &=& X_{\alpha}
\biggl[ \bigl(u_{\alpha}(X) - u(X,X) \bigr)
- \biggl(\eta_{\alpha}^{2} X_{\alpha} - \sum_{\beta} \eta_{\beta}^{2}
X_{\beta}^{2} \biggr) \biggr] \,dt\nonumber\\[-8pt]\\[-8pt]
&&{} + X_{\alpha} \Bigl[\eta_{\alpha} \,dW_{\alpha} - \sum\eta_{\beta}
X_{\beta} \,dW_{\beta} \Bigr].\nonumber
\end{eqnarray}

As it turns out, if a game possesses an interior ESS and the shocks are
mild enough, the solution paths $X(t)$ of the (single-population)
replicator dynamics will be recurrent (Theorem 2.1 in \cite{Im05}).
Theorem \ref{thm:stability} rules out such behavior in the case of
strict equilibria (the multi-population analogue of ESS), but does not
answer the following question: if the underlying game only has mixed
equilibria, will the solution $X(t)$ of the dynamics (\ref{eq:SRD}) be
recurrent?

This question is equivalent to showing that a profile $x$ is
stochastically asymptotically stable in the replicator equations (\ref
{eq:SRD}), (\ref{eq:SLRD}) only if it is a strict equilibrium. Since
Theorem \ref{thm:stability} provides the converse ``if'' part, an
answer in the positive would yield a strong equivalence between
stochastically stable states and strict equilibria; we leave this
direction to be explored in future papers.

\subsection*{\texorpdfstring{It\^{o} vs. Stratonovich}{Ito vs. Stratonovich}}

For comparison purposes (but also for simplicity), let us momentarily
assume that the noise coefficients $\eta_{i\alpha}$ do not depend on
the state $X(t)$ of the game. In that case, it is interesting (and very
instructive) to note that the SDE (\ref{eq:score}) remains unchanged if
we use Stratonovich integrals instead of It\^{o} ones:
%
\begin{equation}
dU_{i\alpha}(t) = u_{i\alpha}(X(t)) \,dt + \eta_{i\alpha}\, \partial
W_{i\alpha}(t).
\end{equation}
Then, after a few calculations, the corresponding replicator equation reads
%
\begin{equation}
\label{eq:Stratonovich}\qquad
\partial X_{i\alpha} = X_{i\alpha} \bigl(u_{i\alpha}(X) - u_{i}(X) \bigr)\,dt
+ X_{i\alpha} \Bigl(\eta_{i\alpha} \,\partial W_{i\alpha} - \sum\eta
_{i\beta}
X_{i\beta}\, \partial W_{i\beta} \Bigr).
\end{equation}

The form of this last equation is remarkably suggestive. First, it
highlights the role of the modified game $\widetilde{u}_{i\alpha} =
u_{i\alpha} + \frac{1}{2}\eta_{i\alpha}^{2}$ even more crisply than
(\ref{eq:SRD}): the payoff terms are completely decoupled from
the noise, in contrast to what one obtains by introducing Stratonovich
perturbations in the evolutionary setting \cite{Im09,KP06}. Secondly,
one can seemingly use this simpler equation to get a much more
transparent proof of Proposition \ref{prop:dominated}: the estimates
for the cross entropy terms $G_{q_{i}-q_{i}'}$ are recovered almost
immediately from the Stratonovich dynamics. However, since (\ref
{eq:Stratonovich}) takes this form only for constant coefficients $\eta
_{i\alpha}$ (the general case is quite a bit uglier), we chose the
route of consistency and employed It\^{o} integrals throughout our paper.

\subsection*{Applications in network design}

Before closing, it is worth pointing out the applicability of the above
approach to networks where the presence of noise or uncertainty has two
general sources. The first of these has to do with the time variability
of the connections which may be due to the fluctuations of the link
quality because of mobility in the wireless case or because of external
factors (e.g., load conditions) in wireline networks. This variability
is usually dependent on the state of the network and was our original
motivation in considering noise coefficients $\eta_{i\alpha}$ that are
functions of the players' strategy profile; incidentally, it was also
our original motivation for considering randomly fluctuating payoffs in
the first place: travel times and delays in traffic models are not
determined solely by the players' choices, but also by the fickle
interference of nature.

The second source stems from errors in the measurement of the payoffs
themselves (e.g., the throughput obtained in a particular link) and
also from the lack of information on the payoff of strategies that were
not employed. The variability of the noise coefficients $\eta_{i\alpha
}$ again allows for a reasonable approximation to this problem. Indeed,
if $\eta_{i\alpha}\dvtx\Delta\to{\mathbb R}$ is continuous and satisfies
$\eta
_{i\alpha}(x_{-i};\alpha) =0$ for all $i\in\mathcal{N}, \alpha\in
\mathcal{S}_{i}$,
this means that there are only errors in estimating the payoffs of
strategies that were not employed (or small errors for pure strategies
that are employed with high probability). Of course, this does not yet
give the full picture [one should consider the discrete-time dynamical
system (\ref{eq:discretescore}) instead where the players' \textit
{actual} choices are considered], but we conjecture that our results
will remain essentially unaltered.

\section*{Acknowledgments}

We would like to extend our gratitude to the anonymous referee for his
insightful comments and to David Leslie from the university of Bristol
for the fruitful discussions on the discrete version of the exponential
learning model.

Some of the results of Section \ref{sec:dominated} were presented in
the conference ``Game Theory for Networks'' in Bo\u{g}azi\c{c}i
University, Istanbul, May 2009 \cite{GameNets09}.


%
\printaddresses

\end{document}